\documentclass[12pt]{iopart}

\usepackage[hidelinks]{hyperref}
\usepackage{iopams}
 \expandafter\let\csname equation*\endcsname\relax
\usepackage[english]{babel}
\usepackage{changes}
\expandafter\let\csname endequation*\endcsname\relax

\usepackage{amsmath}

\newcommand{\argmin}[1]{\underset{#1}{\operatorname{arg}\,\operatorname{min}}\;}

\usepackage{graphicx}
\usepackage{fancyhdr}
\usepackage{booktabs,array,ragged2e}
\usepackage[font=footnotesize]{caption}
\usepackage[font=scriptsize]{subcaption}
\usepackage{appendix}
\usepackage{algorithm,algorithmic}
\usepackage{cite}
\usepackage{wrapfig}
\newcolumntype{C}[1]{>{\centering}m{#1}}
\newcolumntype{L}[1]{>{\raggedright}m{#1}}
\usepackage{url}
\addto\extrasenglish{}
\addto\extrasenglish{}
\usepackage{marginnote}

 \newtheorem{theorem}{Theorem}

\newtheorem{assumption}{Assumption}

\newtheorem{definition}{Definition}

\newtheorem{lemma}{Lemma}

\newtheorem{remark}[theorem]{Remark}

\newenvironment{proof}[1][Proof]{\textit{#1.} }{\ \rule{0.5em}{0.5em}}
\nonstopmode

\begin{document}

\title[~]{Enhancing joint reconstruction and segmentation with non-convex Bregman iteration}

\author{Veronica Corona$^1$, Martin Benning$^2$, Matthias J. Ehrhardt$^3$, Lynn F. Gladden$^4$, Richard Mair$^5$, Andi Reci$^4$, Andrew J. Sederman$^4$, Stefanie Reichelt$^5$, Carola-Bibiane Sch{\"o}nlieb$^1$ }

\address{$^1$ Department of Applied Mathematics and Theoretical Physics, University of Cambridge, United Kingdom\\
$^2$ School of Mathematical Sciences, Queen Mary University of London, United Kingdom\\
$^3$ Institute for Mathematical Innovation, University of Bath, United Kingdom \\
	$^4$ Department of Chemical Engineering and Biotechnology, University of Cambridge, United Kingdom\\
	$^5$Cancer Research UK Cambridge Institute, University of Cambridge, United Kingdom}
\ead{vc324@cam.ac.uk}
\vspace{10pt}

\begin{abstract}
All imaging modalities such as computed tomography (CT), emission tomography and magnetic resonance imaging (MRI) require a reconstruction approach to produce an image. A common image processing task for applications that utilise those modalities is image segmentation, typically performed posterior to the reconstruction. Recently, the idea of tackling both problems jointly has been proposed. We explore a new approach that combines reconstruction and segmentation in a unified framework. 
We derive a variational model that consists of a total variation regularised reconstruction from undersampled measurements and a Chan-Vese based segmentation. We extend the variational regularisation scheme to a Bregman iteration framework to improve the reconstruction and therefore the segmentation. We develop a novel alternating minimisation scheme that solves the non-convex optimisation problem with provable convergence guarantees. Our results for synthetic and real data show that both reconstruction and segmentation are improved compared to the classical  sequential approach.
\end{abstract}

%
%
%
%
%

\section{Introduction}
\label{sec:intro}
\begin{figure}[t]
\centering
\begin{subfigure}{0.25\textwidth}
\centering
 \includegraphics[width=\textwidth]{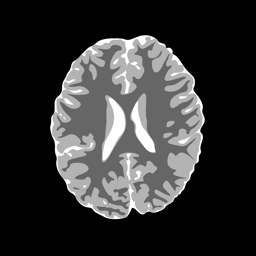}
        \caption{Groundtruth}
        \label{fig:gtbrain}
    \end{subfigure}
~
\begin{subfigure}{0.25\textwidth}
\centering
 \includegraphics[width=\textwidth]{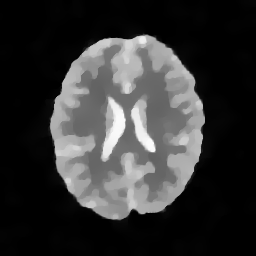}
        \caption{Sequential reconstruction}
        \label{fig:brainTV}
    \end{subfigure}
    ~
    \begin{subfigure}{0.25\textwidth}
    \centering
 \includegraphics[width=\textwidth]{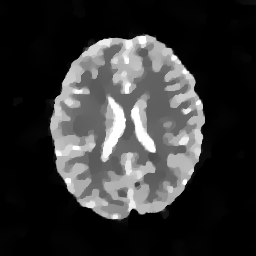}
        \caption{Joint reconstruction}
        \label{fig:brainJoint}
    \end{subfigure}

\begin{subfigure}{0.25\textwidth}
\centering
 \includegraphics[width=\textwidth]{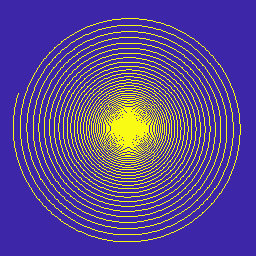}
        \caption{Sampling matrix}
        \label{fig:samp15}
    \end{subfigure}
    ~
            \begin{subfigure}{0.25\textwidth}
        \centering
    \includegraphics[width=\textwidth]{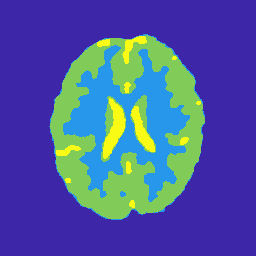}
        \caption{Sequential segmentation}
        \label{fig:brainTVseg}
    \end{subfigure}
    ~
    \begin{subfigure}{0.25\textwidth}
    \centering
 \includegraphics[width=\textwidth]{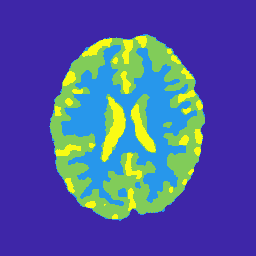}
        \caption{Joint segmentation}
        \label{fig:brainJointseg}
    \end{subfigure} 
    \caption{Sequential approach (left) versus unified approach (right). Combining reconstruction and segmentation in a single unified approach improves both the reconstructed image and its segmentation. See \autoref{fig:brain_whole} for more details.}
    \label{fig:brain}
\end{figure}

Image reconstruction plays a central role in many imaging modalities for medical and non-medical applications.
The majority of imaging techniques deal with incomplete data and noise, making the inverse problem of reconstruction severely ill-posed. Based on compressed sensing (CS) it is possible to tackle this problem by exploiting prior knowledge of the signal \cite{Candes2006,Donoho2006,LustigDonohoPauly2007}. Nevertheless, reconstructions from very noisy and undersampled data will present some errors that will be propagated into further analysis, e.g. image segmentation. Segmentation is an image processing task used to partition the image into meaningful regions. Its goal is to identify objects of interest, based on contours or similarities in the interior. Typically segmentation is performed after reconstruction, hence its result strongly depends on the quality of the reconstruction. Recently the idea of combining reconstruction and segmentation has become more popular. The main motivation is to avoid error propagations that occur in the sequential approach by estimating edges simultaneously from the data, ultimately improving the reconstruction. In this paper, we propose a new model for joint reconstruction and segmentation from undersampled MRI data. The underlying idea is to incorporate prior knowledge about the objects that we want to segment in the reconstruction step, thus introducing additional regularity in our solution. In this unified framework, we expect that the segmentation will also benefit from sharper reconstructions. We demonstrate that our joint approach improves the reconstruction quality and yields better segmentations compared to sequential approaches. In  \autoref{fig:brain}, we consider a brain phantom from which we simulated the undersampled \textit{k}-space data and added Gaussian noise. . 
\autoref{fig:brainTV} and \ref{fig:brainTVseg} 
present reconstructions and segmentations obtained with the sequential approaches, while \autoref{fig:brainJoint} and \ref{fig:brainJointseg} show the results for our joint approach. The reconstruction using our method shows clearly more details and it is able to detect finer structures that are not recovered with the classical separate approach. As a consequence, the joint segmentation is also improved. In the following section we present the mathematical models that we used in our comparison. We investigated the performance of our model for two different applications: bubbly flow and cancer imaging. We show that both reconstruction and segmentation benefit from this method, compared to the traditional sequential approaches, suggesting that error propagation is reduced. 

\paragraph{Our contribution.} In our proposed joint method, we obtain an image reconstruction that preserves its intrinsic structures and edges, possibly enhancing them, thanks to the joint segmentation, and simultaneously we achieve an accurate segmentation. We consider the edge-preserving total variation regularisation for both the reconstruction and segmentation term using Bregman distances. In this unified Bregman iteration framework, we have the advantage of improving the reconstruction by reducing the contrast bias in the TV formulation, which leads to more accurate segmentation. In addition, the segmentation constitutes another prior for the reconstruction by enhancing edges of the regions of interest. Furthermore, we propose a non-convex alternating direction algorithm in a Bregman iteration scheme for which we prove global convergence. \\
\newline
\indent The paper is organised as follows. In \autoref{sec:MRI_RecSeg} we describe the problems of MRI reconstruction and region-based segmentation. We then introduce our joint reconstruction and segmentation approach in a Bregman iteration framework. This section also contains a detailed comparison of other joint models in the literature. In \autoref{sec:optimisation} we study the non-convex optimisation problem and present the convergence analysis for this class of problems. Finally in \autoref{sec:results} we present numerical results for MRI data for different applications. Here we investigate the robustness of our model by testing the undersampling rate up to its limit and by considering different noise levels.  

\section{MRI reconstruction and segmentation}
\label{sec:MRI_RecSeg}
In the following section we introduce the mathematical tools to perform image reconstruction and image segmentation. In this work, we focus on the specific MRI application; however, our proposed joint method can be applied to other imaging problems in which the measured data is connected to the image via a linear and bounded forward operator, cf. \autoref{subsec:rec}. Finally we present our model that combines the two tasks of reconstruction and segmentation in a unified framework.
\subsection{Reconstruction} \label{subsec:rec}
In image reconstruction problems, we have the general setting 
\begin{equation}
f=Au+\eta, 
\label{eq:generalforward}
\end{equation}
where $A:\mathbb{X}\to \mathbb{Y}$ is a bounded and linear operator mapping between two vector-spaces. The measured data $f\in \mathbb{Y}$ is usually corrupted by some noise $\eta$ and often only observed partially. In this formulation we are interested in recovering the image $u$ given the data $f$.\\
In this work, we focus on the application of MRI and we refer to the measurements $f$ as the \textit{k}-space data.
In standard MRI acquisitions, the Fourier coefficients are collected in the \textit{k}-space by radio-frequency (RF) coils. Because the \textit{k}-space data is acquired sequentially, the scanning time is constrained by physical limitations of the imaging system. 
One of the most common ways to perform fast imaging consists of undersampling the \textit{k}-space; this, however, only yields satisfactory results if the dimension of the parameter space can implicitly be reduced, for example by exploiting sparsity in certain domains. In the reconstruction, this assumption is incorporated in the regularisation term. 
Let $\Omega:=\{1, \dots , n_1 \} \times \{1, \dots , n_2 \}$ with $n_1, n_2 \in \mathbb{N}$ be a discrete image domain. Let $f= (f_i)_{i=1}^m \in \mathbb{C}^m$ with $m\ll n=n_1 n_2$ be our given undersampled \textit{k}-space data, where $f_i \in \mathbb{C}$ are the measured Fourier coefficients that fulfil the relationship (\ref{eq:generalforward}) with $A= SF$.
The operator $A$ is now composed by $\mathcal{S}:\mathbb{C}^n \to \mathbb{C}^m$, which is a sampling operator that selects  $m$ measurements from the $Fu$ data according to the locations provided by a binary sampling matrix  (see e.g. \autoref{fig:samp15}), where $F$ is the discrete Fourier transform. 
In MRI, the noise $\eta$ is drawn from a complex-valued Gaussian distribution with zero mean and standard deviation $\sigma$ \cite{Macovski1996}. 
\newline
In problem (\ref{eq:generalforward}) for MRI, the aim is to recover the image $u \in \mathbb{C}^n$ from the data. However, in this work we follow the standard assumption that in many applications we have negligible phase, i.e. we are working with real valued, non-negative images. Therefore, we are only interested in $u \in \mathbb{R}^n$; hence we consider the MRI forward operator as $A : \mathbb{R}^{n} \to \mathbb{C}^m$ and its adjoint $A^*:\mathbb{C}^m \to \mathbb{R}^{n} $ as modelled in \cite{Ehrhardt2016mri}. Problem (\ref{eq:generalforward}) is ill-posed due to noise and incomplete measurements. The easiest approach to approximate (\ref{eq:generalforward}) is to compute the solution, for which the missing entries are replaced with zero
\begin{equation*}
u_z=  A^* f. 
\label{eq:zerofilling}
\end{equation*}
However, images reconstructed with this approach will suffer from aliasing artifacts because undersampling the \textit{k}-space violates the Nyquist--Shannon sampling theorem. 
Therefore, we consider a mathematical model that incorporates prior knowledge by using a variational regularisation approach.  A popular model is to find an approximate solution for $u$ as a minimiser of the Tikhonov-type regularisation approach
\begin{equation}
u^* \in \argmin{u} \Big\{ \frac{1}{2} \| Au - f\|_{2}^{2} + \alpha  J(u) \Big\}, 
\label{eq:tik}
\end{equation}
where the first term is the data fidelity that forces the reconstruction to be close to the measurements and the second term is the regularisation, which imposes some regularity on the solution. The parameter $\alpha > 0$ is a regularisation parameter that balances the two terms in the variational scheme. 
In this setting, different regularisation functionals $J$ can be chosen (see \cite{BenningGladdenHollandEtAl2014} for a survey of variational regularisation approaches). \\ 
Although problems of the form (\ref{eq:tik}) are very effective, they also lead to a systematic loss of contrast \cite{Meyer2001,StrongChan2003,groundstates}. This is typically observed for common choices of the regulariser $J$, i.e., convex functional. To overcome this problem, \cite{Osher2005} proposed an iterative regularisation method based on the generalised Bregman distance \cite{breg,Kiwiel1997}. The Bregman distance with respect to $J$ is defined as
\begin{equation}
D_{J}^{p^k}(u,u^k)= J(u) - J (u^k) - \langle p^k,u-u^k \rangle
\label{bregmandistance}
\end{equation}
with $p^k \in \partial J(u^k) $, where $\partial J (u^k)$ is called sub-differential and it is a generalisation of the classical differential for convex functions. We replace problem (\ref{eq:tik}) with a sequence of minimisation problems
\begin{equation}
u^{k+1} \in \argmin{u} \Big\{ \frac{1}{2} \| Au - f\|_{2}^{2} + \alpha  D_{J}^{p^k}(u,u^k) \Big\}.
\label{eq:rec_breg_J}
\end{equation}
The update on the subgradient can be conveniently computed by the optimality condition of  (\ref{eq:rec_breg_J})
\begin{equation}
p^{k+1} = p^k - \frac{1}{\alpha} A^* (Au^{k+1} - f).
\label{eq:pk}
\end{equation}

In this work, we will focus on one particular choice for $J$, namely the \textit{total variation}. The total variation (TV) regularisation is a well-known edge-preserving approach, first introduced by Rudin, Osher and Fatemi in \cite{ROF} for image denoising. The TV regularisation, i.e., the 1-norm penalty on a discrete finite difference approximation of the two-dimensional gradient $ \nabla : \mathbb{R}^n \to (\mathbb{R}^2)^n$, that is $ \nabla u(i,j) = ( \nabla_1 u(i,j), \nabla_2 u(i,j) )^T$,  
is in the discrete setting
\begin{equation}
J(u)=\operatorname{TV}(u)= \| \nabla u \|_{2,1} = \sum_{(i,j)\in \Omega} \sqrt{| \nabla_1 u(i,j) |^2 + | \nabla_2 u(i,j) |^2} ,
\label{eq:dtv}
\end{equation}
for the isotropic case. \\
\indent We then consider the Bregman iteration scheme in (\ref{eq:rec_breg_J}) for $J(u)=\operatorname{TV}(u)$. 
This approach is usually carried on by initialising the regularisation parameter $\alpha$ with a large value, producing overregularised initial solutions. At every step $k$, finer details are added. A suitable criterion to stop iterations (\ref{eq:rec_breg_J}) and (\ref{eq:pk}) (see \cite{BenningGladdenHollandEtAl2014}), is the Morozov's discrepancy principle \cite{Morozov1966}. The discrepancy principle suggests to choose the smallest $k \in \mathbb{N}$ such that 
$u^{k+1}$ satisfies 
\begin{equation}
\| f - Au^{k+1}\|_2  \leq \sigma \sqrt{m}
\label{eq:discrepancy}
\end{equation}
where $m$ is the number of samples and $\sigma$ is the standard deviation of the noise in the data. Note that using Bregman iterations, the  contrast is improved and in some cases even recovered exactly, compared to the variational regularisation model. In addition, it makes the regularisation parameter choice less challenging. Note that for different choices of $J$ in \eqref{eq:tik}, e.g., the Mumford-Shah/Potts model \cite{Mumford1989,Ambrosio1990,Chambolle1995,Boysen2009,Pock2009}, we do not have loss of contrast, but we deal with a non-convex NP hard problem, algorithmically more challenging. 

\subsection{Segmentation}
Image segmentation refers to the process of automatically dividing the image into meaningful regions. Mathematically, one is interested in finding a partition $\{\Omega_i\}_{i=1}^l $ of the image domain $\Omega$ subject to $\cup_{i=1}^l \Omega_i = \Omega $ and $ \cap_{i=1}^l \Omega_i = \emptyset$.  One way to do this is to use \textit{region-based} segmentation models, which identify regions based on similarities of their pixels. 
The segmentation model we are considering was originally proposed by Chan and Vese in \cite{Chant2001} and it is a particular case of the piecewise-constant Mumford-Shah model \cite{Mumford1989}. Given an image function $u:\Omega \to \mathbb{R}$, the goal is to divide the image domain $\Omega$ in two separated regions $\Omega_1$ and $\Omega_2=\Omega \setminus \Omega_1$ by minimising the following energy function
\begin{equation*}
 \int_{\Omega_1} (u(x) - c_1)^2 \,dx + \int_{\Omega_2} (u(x) - c_2)^2 \,dx + \beta \cdot Length(C) \to \min_{c_1,c_2,C}
\end{equation*}
where $C$ is the desired contour separating $\Omega_1$ and $\Omega_2$, and the constants $c_1$ and $c_2$ represents the average intensity value of $u$ inside $C$ and outside $C$, respectively. The parameter $\beta$ penalises the length of the contour $C$, controlling the scale of the objects in the segmentation. From this formulation we can make two observations: first, the regions $\Omega_1$ and $\Omega \setminus \Omega_1$ can be represented by the characteristic function
\begin{equation*}
v(x)= 
\begin{cases}
      0, & \text{if }  x \in \Omega_1 \cup C \\
      1, & \text{if } x \in \Omega_2,
    \end{cases}
\end{equation*}
second, the perimeter of the contour identified by the the characteristic function corresponds to its total variation, as shown by the Coarea formula \cite{Ambrosio2000}. This leads to the new formulation
\begin{equation*}
 \int_{\Omega} v(x) (u(x) - c_1)^2 dx + (1-v(x)) (u(x) - c_2)^2 \,dx + \beta \operatorname{TV}(v) \to \min_{c_1,c_2, v \in \{0,1\}}.
\end{equation*}
Even assuming fixed constants $c_1$, $c_2$ the problem is non-convex due to the binary constraint. In \cite{Chan2006} the authors proposed to relax the constraint, allowing $v(x)$ to assume values in the interval $[0,1]$. They showed that for fixed constants $c_1$, $c_2$, global minimisers can be obtained by minimising the following energy
\begin{equation}
 \int_{\Omega} v(x) (u(x) - c_1)^2 dx + (1-v(x)) (u(x) - c_2)^2\, dx + \beta \operatorname{TV}(v) \to \min_{v \in [0,1]}
\label{eq:CEN}
\end{equation}
followed by thresholding, setting $\Sigma = \{x : v(x) \geq \mu \} \; \text{for a.e. }\mu\in [0, 1]$. As the problem is convex but not strictly convex, the global minimiser may not be unique. In practice we obtain solutions which are almost binary, hence the choice of $\mu$ is not crucial.  \\ 
Setting 
\begin{equation*}
s(x)=(u(x) - c_1)^2 - (u(x) - c_2)^2
\end{equation*}
the energy \eqref{eq:CEN} can be written in a more general form as
\begin{equation*}
\int_{\Omega}  v(x) s(x)\, dx + \beta \operatorname{TV}(v) \to \min_{v\in  [0,1]}. 
\end{equation*} 
In this paper, we are interested in the extension of the two-class problem to the multi-class formulation \cite{Lellmann2013}. Following the simplex-constrained vector function representation for multiple regions and its convex relaxation proposed in \cite{Lellmann2009}, we obtain as a special case a convex relaxation of the Chan-Vese model for arbitrary number of regions, which reads
\begin{equation}
\int_{\Omega} \sum_{i=1}^{\ell} v_i(x) (c_i - u(x) )^2\, dx + \beta \operatorname{TV}(v) \to \min_{v\in \mathcal{C}},
\label{eq:seg_multi}
\end{equation}
 where  $\mathcal{C}:= \{ v:\Omega \to \mathbb{R}^\ell \,\lvert \, v(x)\geq 0, \sum_{i=1}^\ell v_i(x)=1\}$ is a convex set which restricts $v(x)$ to lie in the standard probability simplex. 
As in the binary case, the constants $c_i$ describe the average intensity value inside region $i$. In this case we consider the vector-valued formulation of TV
 \begin{equation*}
 \operatorname{TV}(v)=\int_{\Omega} \sqrt{\| \nabla v_1 \|^2 + \dots +\| \nabla v_\ell \|^2}\, dx.
 \end{equation*}
\subsection{Joint reconstruction and segmentation }
\label{subsec:joint}
MRI reconstructions from highly undersampled data are subject to errors, even when prior knowledge about the underlying object is incorporated in the mathematical model. 
It is often required to find a trade-off between filtering out the noise and retrieving the intrinsic structures while preserving the  intensity configuration 
and small details. As a consequence, segmentations in the presence of artifacts are likely to fail. \\
\indent In this paper, we propose to solve the two image processing tasks of reconstruction and segmentation in a unified framework. The underlying idea is to inform the reconstruction with prior knowledge of the regions of interest, and simultaneously update this belief according to the actual measurements. Mathematically, given the under-sampled and noisy \textit{k}-space data $f$, we want to recover the image $u \colon \Omega \to \mathbb{R}$ and compute its segmentation $v$ in $\ell$ disjoint regions, by solving the following problem
\begin{equation}
\begin{aligned}
(u,v) &= \argmin{u,v}  \underbrace{\frac{1}{2} \| Au - f\|_{2}^{2} + \alpha \operatorname{TV}(u) }_\text{reconstruction}\\
 &+\underbrace{ \delta \sum_{i=1}^n \sum_{j=1}^{\ell} v_{ij} (c_j - u_i )^2 + \beta \operatorname{TV}(v) + \imath_\mathcal{C}(v)}_\text{segmentation}.
 \end{aligned}
 \label{eq:jointTV}
 \end{equation}
where $\imath_\mathcal{C}(v)$ is the characteristic function over $\mathcal{C}  := \{ v: \mathbb{R}^n \to \mathbb{R}^\ell \,\lvert \, v_{ij}\geq 0, \sum_{j=1}^\ell v_{ij}=1, \forall i  \in \{1,\dots,n\}\}$, and $\alpha$, $\beta$, $\delta >0$ are some regularisation parameters.
However, instead of solving \eqref{eq:jointTV}, we 
 consider the iterative regularisation procedure using Bregman distances. The main motivation is to exploit the contrast enhancement aspect for the reconstruction thanks to the Bregman iterative scheme. By improving the reconstruction, the segmentation is in turn refined. Therefore, we replace \eqref{eq:jointTV} with the following sequence of minimisation problems for $k=0,1,2, \dots$
\begin{subequations}
\begin{align}
u^{k+1} &= \argmin{u} \frac{1}{2} \| Au - f\|_{2}^{2} + \alpha  D_{\operatorname{TV}}^{p^k}(u,u^k) +\delta  \sum_{i=1}^n \sum_{j=1}^{\ell} v_{ij} (c_j - u_i )^2 \label{eq:joint_u}\\
p^{k+1} &= p^k - \frac{1}{\alpha}  \left( A^*(Au^{k+1} - f) - 2 \delta  \sum_{j=1}^\ell v_j^k (u^{k+1}-c_j)  \right) \label{eq:joint_p}\\
v^{k+1} &\in \ \argmin{v} \delta  \sum_{i=1}^n \sum_{j=1}^{\ell} v_{ij} (c_j - u_i^{k+1} )^2 + \chi_\mathcal{C}(v)+ \beta D_{\operatorname{TV}}^{q^k}(v,v^k) \label{eq:joint_v}\\
q^{k+1} &=q^k - \frac{\delta }{\beta} (c_j-u^{k+1})^2.\label{eq:joint_q}  
\end{align}
\label{eq:joint}
\end{subequations} 
\noindent \hspace{-0.167cm}Note that \eqref{eq:joint} solves a problem different from \eqref{eq:jointTV}. Assuming that a minimiser exists, the model \eqref{eq:joint} converges to a minimiser of $$ \frac{1}{2} \| Au - f\|_{2}^{2} + \delta \sum_{i=1}^n \sum_{j=1}^{\ell} v_{ij} (c_j - u_i^{k+1} )^2 \, , $$ 
as we will show in \autoref{sec:convergenceanalysis}. In case of noisy data $f$ this is not desirable, so that we combine the iteration with a stopping criterion in order to form a regularisation method.\\
\indent This model combines the reconstruction approach described in \eqref{eq:rec_breg_J} and the discretised multi-class segmentation in \eqref{eq:seg_multi} with a variation in the regularisation term, which is now embedded in the Bregman iteration scheme. In \cite{Zeune2017} the authors used Bregman distances for the Chan-Vese formulation \eqref{eq:CEN}, combined with spectral analysis, to produce multiscale segmentations. \\ 
\indent As described in the previous subsection, the parameters $\alpha$ and $\beta$ describe the scale of the details in $u$ and the scale of the segmented objects in $v$. By integrating the two regularisations into the same Bregman iteration framework, we obtain that these scales are now determined by the iteration $k+1$. At the first Bregman iteration $k=0$, when $\alpha$ is very large, we obtain an over-smoothed $u^{1}$, and the value of $\beta$ is not very important. Intuitively, $u^1$ is almost piecewise constant with small total variation and a broad range of values of $\beta$ may lead to very similar segmentations $v^1$. However, at every iteration $k+1$, finer scales are added to the solution with the update $p^{k+1}$. Accordingly, with the update $q^{k+1}$, which is independent of $v^{k+1}$, the segmentation keeps up with the scale in the reconstructed image $u^{k+1}$.\\ 
\indent The novelty of this approach is also represented by the role of the parameter $\delta >0$. This parameter weighs the effect of the segmentation in the reconstruction, imposing regularity in $u$ in terms of sharp edges in the regions of interest. In \autoref{sec:results} we show how different ranges of $\delta$ affects the reconstruction (see \autoref{fig:deltas}). Intuitively, large values of $\delta$ force the solution $u$ to be close to the piecewise constant solution described by the constants $c_i$. This is beneficial in applications where MRI is a means to extract shapes and sizes of underlying objects, (e.g. bubbly flow in \autoref{subsec:bubbly}). On the other hand, with very small $\delta$, the segmentation has little impact and the solutions for $u$ are close to the ones obtained by solving the individual problem \eqref{eq:rec_breg_J}. Instead, intermediate values of $\delta$ impose sharper boundaries in the reconstruction while preserving the texture. \\
\indent Obviously, we need to stop the iteration before the residual brings back noise from the data $f$. As we cannot use Morozov discrepancy principle in this case (due to the fact that $ \| Au^k - f \|_2$ will rather increase due to the effect of the coupling term controlled by the parameter $\delta$), we stop when two consecutive iterates in $v$ are smaller than a certain tolerance, $\|v^{k+1}-v^k\| < tol$, following the observation that the rate at which $u^{k+1}$ changes close to the optimal solution is low, in contrary to more abrupt changes  at the beginning of the Bregman iteration and later on when it starts to add noise. \\

\indent Clearly, problem \eqref{eq:joint} is non-convex in the joint argument $(u,v)$ due to the coupling term. However, it is convex in each individual variable. We propose to solve the joint problem by iteratively alternating the minimisation with respect to $u$ and to $v$ (see \autoref{sec:optimisation} for numerical optimisation and convergence analysis). \\
\subsection{Comparison to other joint reconstruction and segmentation approaches} \label{subsec:comparison}
In this section we will provide an overview of some existing simultaneous reconstruction and segmentation (SRS) approaches with respect to different imaging applications. 

\paragraph{CT/SPECT.} Ramlau and Ring \cite{ramlau} first proposed a simultaneous reconstruction and segmentation model for CT, that was later extended to SPECT in \cite{klann} and to limited data tomography \cite{Klann1}. In these work, the authors aim to simultaneously reconstruct and segment the data acquired from SPECT and CT. CT measures the \textit{mass density distribution} $\mu$, that represents the attenuation of x-rays through the material; SPECT measures the \textit{activity distribution} $f$ as the concentration of the radio tracer injected in the material. Given the two measurements $z^{\delta}$
and $y^{\delta}$, from CT and SPECT, they consider the following energy functional
\begin{equation*}
E(f,\mu,\Gamma^f,\Gamma^{\mu}) = \|A(f,\mu) - y^{\delta} \|^2 + \beta \|R \mu -z^{\delta}\|^2 +\alpha ( Length( \Gamma^f)+ Length(\Gamma^{\mu})).
\end{equation*}
 They propose a joint model based on a Mumford-Shah like functional, in which the reconstructions of $\mu$ and $f$ and the given data are embedded in the data term in a least squares sense. The operators $A$ and $R$ are the attenuated Radon transform (SPECT operator) and the Radon transform (CT operator), respectively. The penalty term is considered to be a multiple of the lengths of the contours of $\mu$, $\Gamma^{\mu}$ and the contours of $f$, $\Gamma^{f}$. These boundaries are modelled using level set functions. In these segmented partitions of the domain, $\mu$ and $f$ are assumed to be piecewise constant. 
The optimisation problem is then solved alternatively with respect to the functional variables $f$ and $\mu$ with fixed geometric variables $\Gamma^{\mu}$ and $\Gamma^{f}$ and the other way around. \\
\newline
\indent In \cite{0266-5611-32-10-104002} the simultaneous reconstruction and segmentation is applied to dynamic SPECT imaging, which solves a variational framework consisting of a Kullback-Leibler (KL) data fidelity and different regulariser terms to enforce  sharp edges and sparsity for the segmentation and smoothness for the reconstruction. The cost function is 
\begin{equation*}
E(u,c) = KL\big(R(u\cdot c),g\big) + \alpha \sum_{k=1}^K \|\nabla u_k\| 
 + \beta \sum_{k=1}^K \|u_k\|_1 
+ \frac{\delta}{2} \sum_{k=1}^K \| \frac{\partial}{\partial t} c_k \|^2_2. 
\end{equation*}
Given the data $g$, they want to retreive the concentration curves $c_k(t)$ in time for K disjoint regions and their indication functions $u_k(x)$ in space. The optimisation is carried out alternating the minimisation over $u$ having $c$ fixed and then over $c$ having $u$ fixed. \\
\newline
\indent In \cite{Lauze2017} they propose a variational approach for reconstruction and segmentation of CT images, with limited field of view and occluded geometry. The cost function
\begin{equation*}
\begin{aligned}
E(u,c,v)=\frac{1}{2} \|Ax-y\|^2 + \alpha \| \nabla u\| + \frac{\beta}{2} \left( \lambda \sum_i^n \sum_{k=1}^K v_{ik} (u_i-c_k)^2 + \frac{1}{2} \|Dv\|^2_2 \right)\\
\end{aligned}
\end{equation*}
s.t. a box constraint on the image values $x$ and the simplex constraint on the labelling function $v$. The operator $A$ is the undersampled Radon transform modelling the occluded geometry and $y$ is the given data. The second term is the edge-preserving regularisation term for $u$, the third term is the segmentation term which aims at finding regions in $u$ that are  close to the value $c_k$ in region $k$. The operator $D$ is the finite difference approximation of the gradient. The non-convex problem is solved by alternating minimisation between updates of $u,v,c$. 
\paragraph{PET and Transmission Tomography.} In \cite{van}, the authors propose a maximum likelihood reconstruction and doubly stochastic segmentation for emission and transmission tomography. In their model they use a Hidden Markov Measure Field Model (HMMFM) to estimate the different classes of objects from the given data $r$. They want to maximise the following cost function
\begin{equation*}
E(u,p,\theta)= \log P(r|u) + \log P(u|p,\theta) + \log P(p).
\end{equation*}
The first term is the data likelihood which will be modelled differently for emission and transmission tomography. The second term is the conditional probability or class fitting term, for which they use HMMFM. The third term is the regularisation on the HMMFM. The optimisation is carried out in three steps, where first they solve for $u$ (image update) fixing $p, \theta$, then for $p$, holding $u,\theta$ (measure field update) and finally for $\theta$ (parameter update) having $u,p$ fixed. \\
\newline
\indent A variant of this method has been presented in \cite{yiqiu}, in which they incorporate prior information about the segmentation classes through a HMMFM. Here, the reconstruction is the minimisation over a constrained Bayesian formulation that involves a \textit{data fidelity term} as a classical least squares fitting term, a \textit{class fitting term} as a Gaussian mixture for each pixel given $K$ classes and dependent of the class probabilities defined by the HMMFM, and a \textit{regulariser} also dependent of the class probabilities. The model to minimise is
\begin{equation*}
\begin{aligned}
E(u,\delta) =& \lambda_{noise} \|Au-b\|^2_2  -  \sum_{j=1}^N \log \left(  \sum_{k=1}^K \frac{\delta_{jk}}{\sqrt{2 \pi} \sigma_k} \exp \left( - \frac{(u_j -\mu_k)^2}{2 \sigma_k^2} \right) \right) + \lambda_{class} \sum_{k=1}^K R(\delta_k) \\
\text{s.t.} \quad &  \sum_{k=1}^K  \delta_{jk}=1, \quad \delta_{jk} \geq 0,\quad j=1,\dots,N, \quad k=1,\dots, K.
\end{aligned}
\end{equation*}
The operator $A$ will be modelled as the Radon transform in case of CT and $b$ represents the measured data; $N$ is the number of pixel in the image; $\lambda_{noise}$ and $\lambda_{class}$ are the regularisation parameters; $\mu_k, \sigma_k$ are the class parameters. The cost function is non convex and they solve the problem in an alternating scheme where they either update the pixel values or the class probabilities for each pixel.\\
\newline
\indent Storath and others \cite{storath} model the joint reconstruction and segmentation using the Potts Model with application to PET imaging and CT. They consider the variational formulation of the Potts model for the reconstruction. Since the solution is piecewise constant, this directly induces a partition of the image domain, thus a segmentation. Given the data $f$ and an operator $A$ (e.g. Radon transform), the energy functional is in the following form
\begin{equation*}
E(u)= \lambda \| \nabla u \|_0 + \|Au-f\|^2_2
\end{equation*}
where the first term is the jump penalty enforcing piecewise constant solutions and the second term is the data fidelity. As the Potts model is NP hard, they propose a discretisation scheme that allows to split the Potts problem into subproblems that can be solved efficiently and exactly. 
\paragraph{MRI.} In \cite{caballero}, the authors proposed a joint model with application to MRI. Their reconstruction-segmentation model consists of a fitting term and a patch-based dictionary to sparsely represent the image, and a term that models the segmentation as a mixture of Gaussian distributions with mean, standard deviation and mixture weights $\mu$, $\sigma$, $\pi$. Their model is 
\begin{equation*}
E(u,\Gamma,\mu,\sigma,\pi) = \|Au-y\|^2+\lambda \sum_{n=1}^N \|R_n u - D \gamma_n\|^2 - \beta \text{ln} P(u|\mu,\sigma,\pi) \quad \text{s.t.} \; \|\gamma_n\|_0 \leq T \; ~ \forall n
\end{equation*}
where $A$ is the undersampled Fourier transform, $y$ is the given data, $R_n$ is a patch extraction operator, $\lambda$ is a weighting parameter, $T$ is the sparsity threshold, and $\gamma_n$ is the sparse representation of patch $R_n u$ organised as column $n$ of the matrix $\Gamma$. The problem is highly non-convex and it is solved iteratively using conjugate gradient  on $u$, orthogonal matching pursuit on $\Gamma$ and Expectation-Maximisation algorithm on $(\mu,\sigma,\pi)$.
\paragraph{Summary.} Recently, the idea to solve the problems of reconstruction and segmentation simultaneously has become more popular. The majority of these joint methods have been proposed for CT, SPECT and PET data. Mainly they differ in the way they encode prior information in terms of regularisers and how they link the reconstruction and segmentation in the coupling term. Some imposes smoothness in the reconstruction \cite{0266-5611-32-10-104002}, others sparsity in the gradient \cite{ramlau,Lauze2017,storath}, other consider a patch-dictionary sparsifying approach \cite{caballero}. In \cite{storath} they do not explicitly obtain a segmentation, but they force the reconstruction to be piecewise constant. Depending on the application, the coupling term is the data fitting term itself (e.g. SPECT), or  the segmentation term. In \cite{van,yiqiu,caballero} the authors model the segmentation as a mixture of Gaussian distribution, while \cite{Lauze2017} has a a region-based segmentation approach similar to what we propose. However, \cite{Lauze2017} penalises the squared 2-norm of segmentation, imposing spatial smoothness. \\
In our proposed joint approach, we perform reconstruction and segmentation in a unified Bregman iteration scheme, exploiting the advantage of improving the reconstruction, which results in a more accurate segmentation. Furthermore, the segmentation constitutes another prior imposing regularity in the reconstruction in terms of sharp edges in the regions of interest. We propose a novel numerical optimisation problem in a non-convex Bregman iteration framework for which we present a rigorous   convergence result in the following section.
\section{Optimisation}
\label{sec:optimisation}
The cost function \eqref{eq:joint} is non-convex in the joint argument $(u,v)$, but it is convex in each individual variable. To solve this problem we derive a splitting approach where we solve the two minimisation problems in an alternating fashion with respect to $u$ and $v$. We present the general algorithm and its convergence analysis in the next subsection. First, we describe the solution of each subproblem.
\paragraph{Problem in $u$.} 
The problem in $u$ reads 
\begin{equation*}
\begin{aligned}
u^{k+1}= \argmin{u} & \frac{1}{2} \| Au - f\|_{2}^{2} + \alpha  (\operatorname{TV}(u) - \langle p^k,u \rangle) + \delta\sum_{i=1}^n \sum_{j=1}^{\ell} v_{ij}^k (c_j - u_i )^2.
\end{aligned}
\label{eq:problemU}
\end{equation*}
We solve the optimisation for $u$, fixing $v$, using the primal-dual algorithm proposed in \cite{pd,ChambollePock2016,EsserZhangChan2010,PockCremersBischofEtAl2009}. 
We write $F(u)= \|  u \|_1$, $K(u)=\nabla u$ and $G(u) =  \frac{1}{2} \| Au - f\|_{2}^{2} - \alpha  \langle p^k,u \rangle + \delta \sum_{i=1}^n \sum_{j=1}^{\ell} v_{ij}^k (c_j - u_i )^2$ 
and obtain the following iterates for $\theta=1$ and step sizes $\sigma=\tau=0.99 / \| \nabla \|$

\begin{equation*}
\begin{aligned}
y^{n+1} &= \frac {y^n + \sigma \nabla \bar{u}^n}{\max (1, \| y^n + \sigma \nabla \bar{u}^n \| ) } \\
u^{n+1} &=   \frac{  u^n + \tau \nabla \cdot y^{n+1} + 2 \tau \delta \sum_{j=1}^\ell v_j ^k c_j + \tau \alpha p^k +\tau A^* f}{1+ 2 \tau \delta + \tau A^*A}   \\
\bar{u}^{n+1} &= 2u^{n+1}-u^n.
\end{aligned}
\end{equation*}
After sufficiently many iterations we set $u^{k+1}=u^{n+1}$ and compute the update $p^{k+1}$ from the optimality condition of \eqref{eq:problemU} as \eqref{eq:joint_p}.

\paragraph{Problem in $v$.}  The problem in $v$ reads
\begin{equation*}
v^{k+1} = \argmin{v\in C}    \langle v, \delta g - \beta q^k \rangle + \beta \operatorname{TV}(v) 
\end{equation*}
with $ g = \big( (c_1 - u^{k+1})^2 , \dots, (c_\ell- u^{k+1})^2 \big)^T$. 
We now solve a variant of the primal-dual method \cite{pd} as suggested in \cite{PockCremersBischofEtAl2009,ChambolleCremersPock2012}. They consider the general problem including pointwise linear terms of the form
\begin{equation*}
\min_{x\in C} \max_{y\in B} \langle Kx, y \rangle + \langle g,x \rangle - \langle h, y\rangle
\end{equation*}
where $ C \subseteq X$, $ B\subseteq Y$ are closed, convex sets.\\
Setting $K=\nabla$ and $h=0$, $\theta=1$ and step sizes $\sigma=\tau=0.99 / \| \nabla \|$, the updates are 
\begin{equation*}
\begin{aligned}
w^{n+1} &=  \Pi_B \big(w^n+\sigma (\nabla \bar{v}^n-h)\big) \\ 
v^{n+1} &=  \Pi_C \big(v^n + \tau \nabla \cdot (v^{n+1} - \delta g + \beta q^k)\big)\\
\bar{v}^{n+1} &= 2v^{n+1}-v^n.
\end{aligned}
\end{equation*}
At the end, we set $v^{k+1}=v^{n+1}$ and obtain the update $q^{k+1}$ as \eqref{eq:joint_q}.

\subsection{Convergence analysis}
\label{sec:convergenceanalysis}
The proposed joint approach \eqref{eq:joint} is an optimisation problem of the form
\begin{equation}
\min_{u,v} E(u,v) + D_{J_1}^{p^k}(u,u^k)+  D_{J_2}^{q^k}(v,v^k)
\label{eq:P}
\end{equation}
in the general Bregman distance framework for (nonconvex) functions $E: \mathbb{R}^n \times \mathbb{R}^m \to \mathbb{R} \cup \{\infty\}$, for $k \in \{0, \dots, N\}$ and some positive parameters $\alpha$ and $\beta$. The functions $J_1:\mathbb{R}^n \to \mathbb{R} \cup \{ \infty \}$ and $J_2:\mathbb{R}^m \to \mathbb{R} \cup \{ \infty \}$ impose some regularity in the solution. In this work we consider a finite dimensional setting and we refer to the next section for the required definitons. 
To prove global convergence of \eqref{eq:P}, we consider functions that satisfy the Kurdika-\L{}ojasiewicz property, defined below, and we make the following assumptions.
{\definition{(Kurdyka-\L{}ojasiewicz (KL) property).}\label{KL}} Let $F: \mathbb{R}^d \to \mathbb{R}$ be a proper and lower semicontinuous function. 
\begin{itemize}
\item Then the function $F$ is said to have the KL property at $\bar{u}\in \text{dom}(\partial F) := \{ u \in \mathbb{R}^d | \partial F \neq \emptyset \} $ if there exists a constant  $\eta \in (0, \infty ]$, a neighbourhood $N$ of $\bar{u}$ and a concave function $\varphi: [0,\eta) \to \mathbb{R}_{>0}$ that is continuous at $0$ and satisfies $\varphi(0)=0$, $\varphi \in C^1(]0,\eta[)$,  and $\varphi'(s)>0$ for all $s\in ]0,\eta [$, such that for all $u\in N \cap \{ u \in \mathbb{R}^d | F(\bar{u}) < F(u)<F(\bar{u})+\eta\}$ the inequality
\begin{equation*}
\varphi'( F(u) -F(\bar{u}) ) \text{dist}(0,\partial F(u)) \geq 1
\label{eq:KL} \tag{KL}
\end{equation*}
holds.
\item If $F$ satisfies the KL property at each point of $\text{dom}(\partial F)$, $F$ is called a KL function.
\end{itemize}

{\lemma{  ~} } The function $E(u,v)=\frac{1}{2} \| Au - f\|_{2}^{2} + \delta \sum_{i=1}^n \sum_{j=1}^{\ell} v_{ij} (c_j - u_i )^2 $ in our joint problem \eqref{eq:joint} satisfies the KL property over $\mathbb{R}^n \times \mathbb{R}^m$. \\
\begin{proof}
It has been proved in \cite{Lojasiewicz1963} that real-analytic functions satisfy the KL property. The function $E(u,v)$ is polynomial and therefore it is a real-analytic function. \end{proof}
{\assumption \qquad \\  \label{as:1} } 
\vspace{-0.5cm}
\begin{enumerate}
\item $E$ is a $C^1$ function
\item $E>-\infty$ 
\item $E$ is a KL function 
\item $J_i: \mathbb{R}^n \to \mathbb{R}$, $i=1,2$, are proper, lower semi-continuous (l.s.c.)  and strongly convex
\item $J_i$, $i=1,2$,  are KL function
\item for any fixed $v$, the function $u \to E(u,v)$ is convex.  Likewise for any fixed $u$, the function  $v \to E(u,v)$ is convex. 
\item for any fixed $v$, the function $u \to E(u,v)$ is $C^{1}_{L_1(v)}$,  hence the partial gradient is $L_1(v)$-Lipschitz continuous
\begin{equation*}
\| \nabla_u E(u_1,v) - \nabla_u E(u_2,v) \| \leq L_1(v) \| u_1 - u_2 \| \qquad \forall u_1,u_2 \in \mathbb{R}^n.
\end{equation*}
 Likewise for any fixed $u$, the function  $v \to E(u,v)$ is $C^{1}_{L_2(u)}$.
\end{enumerate}
We want to study the convergence properties of the alternating scheme

\begin{subequations}
\begin{align}
u^{k+1} &= \argmin{u} \left\{E(u, v^k) + D_{J_1}^{p^k} (u,u^k) \right\} \label{eq:A1.1}\\
p^{k+1} &= p^k - \nabla_u E(u^{k+1}, v^k) \label{eq:A1.2}\\
v^{k+1} &= \argmin{v} \left\{ E(u^{k+1}, v) + D_{J_2}^{q^k} (v,v^k) \right\}\label{eq:A1.3}\\ 
q^{k+1} &= q^k - \nabla_v E(u^{k+1}, v^{k+1}) \label{eq:A1.4} 
\end{align}
\label{eq:A} 
\end{subequations}

\noindent for initial values $(u^0,v^0)$, $p^0 \in \partial J_1(u^0)$ and $q^0 \in \partial J_2(v^0)$.
\\
We want to show that the whole sequence generated by \eqref{eq:A} converges to a critical point of $E$. \\

\begin{algorithm}[t]
\begin{algorithmic}
\STATE{ \textbf{Initialization:}} $(u^0,v^0)$, $p^0 \in \partial J_1(u^0)$,  $q^0 \in \partial J_2(v^0)$, $N\in \mathbb{N}$
  \FOR {$k=0,1,\dots, N$}
  \STATE $u^{k+1} = \argmin{u} \big\{E(u, v^k) + D_{J_1}^{p^k} (u,u^k) \big\}$
  \STATE $p^{k+1} = p^k - \nabla_u E(u^{k+1}, v^k) $
  \STATE $v^{k+1} = \argmin{v} \big\{ E(u^{k+1}, v) + D_{J_2}^{q^k} (v,v^k) \big\} $
 \STATE $q^{k+1} = q^k - \nabla_v E(u^{k+1}, v^{k+1}) $
  \ENDFOR
\end{algorithmic}
\caption{Alternating splitting method with Bregman iterations for two blocks.}
\label{alg1}
\end{algorithm}
In order for the updates \eqref{eq:A1.1} and \eqref{eq:A1.3} to exist, we want $J$ to be of the form $ J=R+\varepsilon G$ (e.g. $ R = \| \nabla u \|_1$ and $ G= \| u \|^2_2$, see \cite{BenningBetckeEhrhardtEtAl2017})  where $R$ and $G$ fulfil the following assumptions. In practice, we verify that $G$ does not significantly change the reconstruction and segmentation performance for the examples we consider in the next section, for sufficiently small parameter (e.g. $\varepsilon=10^{-3}$). Therefore, in our model \eqref{eq:joint} and in the numerical results we omit it. \\ 
{\assumption{  \qquad \\   } \label{as:2} }
\vspace{-0.5cm}
\begin{enumerate}
\item The functions $G_1:\mathbb{R}^n \to \mathbb{R}$ and $G_2:\mathbb{R}^m \to \mathbb{R}$ are strongly convex with constants $\gamma_1$ and $\gamma_2$, respectively. They have Lipschitz continuous gradient $\nabla G_1$ and $\nabla G_2$ with Lipschitz constant $\delta_1$ and $\delta_2$, respectively.
\item The functions $R_1:\mathbb{R}^n \to \mathbb{R}$ and $R_2:\mathbb{R}^m \to \mathbb{R}$ are proper, l.s.c.\ and convex.
\end{enumerate}

For $J_i=\alpha_i R_i+\varepsilon_i G_i$,  $i \in \{1,2\}$, we can write \eqref{eq:A} as
\begin{subequations}
\begin{align}
u^{k+1} &= \argmin{u} \big\{E(u, v^k) + \alpha_1 D_{R_1}^{p^k} (u,u^k)+ \varepsilon_1 D_{G_1}(u,u^k) \big\} \label{eq:A2.1}\\
p^{k+1} &= p^k - \frac{1}{\alpha_1} \Big(\nabla_u E(u^{k+1}, v^k) + \varepsilon_1 \big( \nabla G_1(u^{k+1} )- \nabla G_1(u^k) \big) \Big)  \label{eq:A2.2}\\ 
v^{k+1} &= \argmin{v} \big\{ E(u^{k+1}, v) + \alpha_2 D_{R_2}^{q^k} (v,v^k)+ \varepsilon_2 D_{G_2}(v,v^k)\big\}\label{eq:A2.3} \\
q^{k+1} &= q^k - \frac{1}{\alpha_2} \Big(\nabla_v E(u^{k+1}, v^{k+1}) +  \varepsilon_2 \big( \nabla G_2(v^{k+1} )- \nabla G_2(v^k)\big) \Big). \label{eq:A2.4}
\end{align}
\label{eq:A2} 
\end{subequations}
{\theorem{(Global convergence). \label{theorem}}} Suppose $E$ is a KL function for any $z^k=(u^k,v^k) \in \mathbb{R}^n \times \mathbb{R}^m$ and $r^k=(p^k,q^k)$ with $p^k \in \partial R_1(u^k)$, $q^k \in \partial R_2(v^k)$. Assume Assumptions \autoref{as:1} and \ref{as:2} hold. Let $\{ z^k \}_{k\in \mathbb{N}}$ and $\{r^k \}_{k\in \mathbb{N}}$ be sequences generated by \eqref{eq:A2}, which are assumed to be bounded. Then
\begin{enumerate}
\item The sequence $\{z^k\}_{k \in \mathbb{N}}$ has finite length, that is
\begin{equation}
\sum_{k=0}^{\infty} \|  z^{k+1} - z^k \| < \infty.
\end{equation}
\item The sequence $\{ z^k\}_{k\in \mathbb{N}}$ 
  converges to a critical point $\bar{z}$ of $E$.
\end{enumerate}

\subsection{Proof of \autoref{theorem}}

In the following we are going to show global convergence of this algorithm. The first step in our convergence analysis is to show a sufficient decrease property of a surrogate of the energy function \eqref{eq:P} and a subgradient bound of the norm of the iterates gap. 
We first recall the following definitions.
{\definition{(Convex Conjugate). \label{def:cc}}} Let $G$ be a proper, l.s.c.\ and convex function. Then its convex conjugate $G^*:\mathbb{R}^n \to \mathbb{R} \cup \{\infty \}$ is defined as 
\begin{equation*}
G^*(p):= \sup_{u\in \mathbb{R}^n} \{\langle u,p\rangle -G(u) \},
\end{equation*}
for all $p \in \mathbb{R}^n$.
{\lemma{~ \label{subgradient}}}Let $G$ be a proper, l.s.c.\ and convex function and $G^*$ its convex conjugate. Then for all arguments $u\in \mathbb{R}^n$ with corresponding subgradients $p\in \partial G(u)$ we know
\begin{itemize}
\item $\langle u,p \rangle = G(u) +G^*(p)$,
\item $p\in \partial G(u)$ is equivalent to $u\in \partial G^*(p)$.
\end{itemize}
From Lemma \autoref{subgradient} we can rewrite the Bregman distance in \eqref{bregmandistance} as follows
\begin{equation}
D_J^{p^k}(u,u^k)=J(u)+J^*(p^k)-\langle u, p^k \rangle,
\end{equation}
where we can see that now it does not depend on $u^k$ anymore, but it can be defined as a function of $u$ and $p^k$ only, $D_J(u,p^k)$.

{\definition{(Strong convexity).}} Let $G$ be a proper, l.s.c.\ and convex function. Then $G$ is said to be $\gamma$-strongly convex if there exists a constant $\gamma$ such that
\begin{equation*}
D_G^p(u,v) \geq \frac{\gamma}{2} \|u -v \|^2
\label{eq:strong_convex}
\end{equation*}
holds true for all $u,v \in \operatorname{dom}(G)$ and $q \in \partial G(v)$.

{\definition{(Symmetric Bregman distance).}} Let $G$ be a proper, l.s.c.\ and convex function. Then the symmetric generalised Bregman distance $D_G^{symm} (u,v)$ is defined as
\begin{equation*}
D_G^{symm} (u,v):= D_G^p(u,v) + D_G^q (v,u) = \langle p-q,u-v \rangle 
\end{equation*}
for $u,v \in \operatorname{dom}(G)$ with $p\in \partial G(u)$ and $q \in \partial G(v)$.
We also observe that in case $G$ is $\gamma$-strongly convex we have
\begin{equation*}
D_G^{symm}(u,v) \geq \gamma\|u -v \|^2.
\label{eq:strong_convex_symm}
\end{equation*}
{\definition{(Lipschitz continuity).}} A function $G: \mathbb{R}^n\to \mathbb{R}$ is (globally) Lipschitz-continuous if there exists a constant $L > 0$ such that
\begin{equation*}
\|G(u)-G(v)\| \leq L \| u -v\|
\label{eq:lipschitz}
\end{equation*}
is satisfied for all $u,v \in \mathbb{R}^n$. \\
\newline
\indent Before we show global convergence, we first define the surrogate functions. 
{\definition{(Surrogate objective).} } Let 
$E, R_i, G_i, i\in\{1,2\}$ satisfy Assumption 1 and Assumption 2, respectively. For any $(u^k,v^k) \in \mathbb{R}^n \times \mathbb{R}^m$ and subgradients $p^k \in \partial R_1(u^k)$ and $q^k \in \partial R_2(v^k)$, we define the following \textit{surrogate objectives} $F$, $F_1$ and $F_2$ 
\begin{equation}
\begin{aligned}
F(u^{k+1},v^{k+1},p^k, q^k) =E(u^{k+1},v^{k+1}) &+ \alpha_1\underbrace{\big(R_1(u^{k+1})+R_1^*(p^k)-\langle u^{k+1},p^k \rangle\big)}_{\text{$= D_{R_1}^{p^{k}}(u^{k+1},u^{k})$} }\\
&+ \alpha_2 \underbrace{\big(R_2(v^{k+1})+R_2^*(q^k)-\langle v^{k+1},q^k \rangle\big)}_{\text{$= D_{R_2}^{q^{k}}(v^{k+1},v^{k})$} },
\end{aligned}
\label{eq:surr}
\end{equation}
\begin{equation}
F_1(u^{k+1},p^k) =E(u^{k+1},v^{k+1}) + \alpha_1\big( R_1(u^{k+1})+R_1^*(p^k)-\langle u^{k+1},p^k \rangle \big), 
\end{equation}
\begin{equation}
 F_2(v^{k+1},q^k) =E(u^{k+1},v^{k+1}) + \alpha_2 \big(R_2(v^{k+1})+R_2^*(q^k)-\langle v^{k+1},q^k \rangle \big).
  \end{equation}
For convenience we will use the following notations
\begin{equation*}
\begin{aligned}
z^k:=&(u^k,v^k) \qquad \forall k \geq 0 \\
r^k:=&(p^k,q^k) \qquad p^k \in \partial R_1(u^k), \: q^k \in \partial R_2(v^k).
\end{aligned}
\end{equation*}
The surrogate function $F$ will then read
\begin{equation*}
F(z^{k+1},r^k)=F(u^{k+1},v^{k+1},p^k, q^k).
\end{equation*}

We can now show the sufficient decrease property of \eqref{eq:surr} for subsequent iterates.
{\lemma{(Sufficient decrease property). \label{lemmaDecrease}}}
The iterates generated by \eqref{eq:A2} satisfy the descent estimate

\begin{equation}
\begin{aligned}
F(z^{k+1},r^k) + \rho_2 \|z^{k+1}-z^k\|^2  
\leq F(z^{k},r^{k-1})
\end{aligned}
\label{eq:decrease}
\end{equation}
In addition we observe 
\begin{equation*}
\begin{aligned}
\lim_{k \to \infty} D_{R_1}^{symm} (u^{k+1},u^k) = 0 \qquad \lim_{k \to \infty}   D_{R_2}^{symm}(v^{k+1},v^k)=0 \\
\lim_{k \to \infty} D_{G_1}^{symm} (u^{k+1},u^k) =0 \qquad \lim_{k \to \infty} D_{G_2}^{symm}(v^{k+1},v^k)=0.
\end{aligned}
\end{equation*}
\begin{proof}
From \eqref{eq:P} we consider the following step for $J_1=\alpha_1R_1+\varepsilon_1G_1$
\begin{equation*}
\begin{aligned}
u^{k+1} =& \argmin{u} \big\{E(u, v^k) + \alpha_1 D_{R_1}^{p^k} (u,u^k) +\varepsilon_1 D_{G_1}(u,u^k) \big\}
\\
=& \argmin{u} \big\{E(u, v^k) + \alpha_1 R(u) + \varepsilon_1 G(u) - \langle  \alpha_1 p^k +\varepsilon_1 \nabla G(u^k) , u -u^k\rangle \big\}.
\end{aligned}
\end{equation*}
Computing the optimality condition we obtain
\begin{equation*}
\alpha_1 (p^{k+1} - p^k) + \nabla_u E (u^{k+1},v^k) + \varepsilon_1 (\nabla G(u^{k+1})- \nabla G(u^k))=0
\end{equation*}
Taking the dual product with $u^{k+1} - u^k$ yields
\begin{equation*}
\begin{aligned}
\alpha_1 \underbrace{ \langle p^{k+1} - p^k,u^{k+1} - u^k \rangle}_{\text{$=D_{R_1}^{symm}(u^{k+1},u^k)$} }
+ \underbrace{ \langle \nabla_u E (u^{k+1},v^k),u^{k+1} - u^k \rangle}_{\text{ $\geq E(u^{k+1}, v^k) -E(u^{k}, v^k)$}} \\
 + \varepsilon_1 \underbrace{ \langle \nabla G_1(u^{k+1})- \nabla G_1(u^k),u^{k+1} - u^k \rangle }_{\text{$=D_{G_1}^{symm}(u^{k+1},u^k)$} }=0.
 \end{aligned}
\end{equation*}
Using the convexity estimate $E(u^{k+1},v^k) - E(u^k,v^k) \leq - \langle \nabla_u E (u^{k+1},v^k),u^{k+1} - u^k \rangle$ we obtain the inequality
\begin{equation*}
\alpha_1 D_{R_1}^{symm}(u^{k+1},u^k) + \varepsilon_1 D_{G_1}^{symm}(u^{k+1},u^k) + E(u^{k+1}, v^k) -E(u^{k}, v^k) \leq 0 
\end{equation*}

\begin{equation*}
\begin{aligned}
\alpha_1 \big(D_{R_1}^{p^k}(u^{k+1},u^k) +D_{R_1}^{p^{k+1}}(u^{k},u^{k+1})\big) + \varepsilon_1 D_{G_1}^{symm}(u^{k+1},u^k)  + E(u^{k+1},v^k) 
\\ 
\leq E(u^{k},v^k).
\end{aligned}
\end{equation*}

\noindent Adding $\alpha_1 D_{R_1}^{p^{k-1}}(u^{k},u^{k-1})$ to both sides, using the strong convexity of $G_1$ and the surrogate function notation, we get

\begin{equation*}
\begin{aligned}
F_1(u^{k+1},p^k)+ \alpha_1 \big(D_{R_1}^{p^{k+1}}(u^{k},u^{k+1}) + D_{R_1}^{p^{k-1}}(u^{k},u^{k-1})\big) + \varepsilon_1 \gamma_1 \|u^{k+1} -u^k\|^2  
\leq F_1(u^{k},p^{k-1}).
\end{aligned}
\end{equation*}
Using the trivial estimate for the Bregman distances, we get the decrease property
\begin{equation*}
\begin{aligned}
F_1(u^{k+1},p^k)+ \varepsilon_1 \gamma_1 \|u^{k+1} -u^k\|^2  
\leq F_1(u^{k},p^{k-1}).
\end{aligned}
\end{equation*}
Similarly for $v$, we obtain

\begin{equation*}
\begin{aligned}
F_2(v^{k+1},q^{k})+ \varepsilon_2 \gamma_2 \|v^{k+1} -v^k \|^2 
\leq F_2(v^{k},q^{k-1}).
\end{aligned}
\end{equation*}
Summing up these estimates, we verify the sufficient decrease property \eqref{eq:decrease}, with positive \\$\rho_2 = \max \{ \varepsilon_1 \gamma_1, \varepsilon_2 \gamma_2  \}$. We also observe 
\begin{equation*}
\begin{aligned}
0 &\leq \Delta^k \leq E(z^k) -E(z^{k+1}).
\end{aligned}
\end{equation*}
with 
\begin{equation*}
\Delta^k := \alpha_1 D_{R_1}^{symm} (u^{k+1},u^k)  + \alpha_2 D_{R_2}^{symm}(v^{k+1},v^k) + \varepsilon_1 D_{G_1}^{symm} (u^{k+1},u^k)+ \varepsilon_2 D_{G_2}^{symm}(v^{k+1},v^k).
\end{equation*}
Summing over $k=0,\dots,N$
\begin{equation*}
\begin{aligned}
\sum_{k=0}^{N} &\Delta^k \leq  \sum_{k=0}^{N}  E(z^k) -E(z^{k+1}) = E(z^0) - E(z^{N+1}) \leq E(z^0) - \inf_z E(z) < \infty.
\end{aligned}
\end{equation*}
Taking the limit $N \to \infty$ implies
\begin{equation*}
\sum_{k=0}^{\infty}  \Delta^k < \infty
\end{equation*}
thus $\lim_{k \to \infty}  D_{R_1}^{symm} (u^{k+1},u^k) = 0$ , $\lim_{k \to \infty} D_{G_1}^{symm}=0$,
$\lim_{k \to \infty}   D_{R_2}^{symm}(v^{k+1},v^k)=0$, \\
$\lim_{k \to \infty}D_{G_2}^{symm}(v^{k+1},v^k)=0$, due to $\alpha_1$, $\alpha_2$, $\varepsilon_1$, $\varepsilon_2 >0$.
\end{proof}
\newline
\indent In order to show that the sequences generated by \eqref{eq:A2} approach the set of critical point we first estimate a bound for the subgradients of the surrogate functions and verify some properties of the limit point set. We first write the subdifferential of the surrogate function as 
\begin{equation}
w^{k+1}:=\begin{pmatrix}  \nabla_u E (u^{k+1},v^{k+1}) + \alpha_1 (p^{k+1}-p^k) \\
 \nabla_v E (u^{k+1},v^{k+1})+ \alpha_2 (q^{k+1}-q^k) \\
 u^k-u^{k+1}\\
 v^k-v^{k+1}
 \end{pmatrix} \in \partial F(z^{k+1},r^k) 
 \label{eq:w}
\end{equation}
with $p^k\in\partial R_1(u^k)$ and $q^k \in \partial R_2(v^k)$ being equivalent to $u^k \in \partial R_1^* (p^k)$ and $v^k \in \partial R_2^* (q^k)$, respectively.
{\lemma{(A subgradient lower bound for the iterates gap). \label{lemmaSubgrBound}}}
Suppose Assumptions \autoref{as:1} and \ref{as:2}  hold. Then the iterates \eqref{eq:A2} satisfy
\begin{equation}
\|w^{k+1}\| \leq \rho_1 \| z^{k+1} - z^k \| 
\label{eq:bound}
\end{equation}

$w^{k+1}\in \partial F(z^{k+1},r^k)$ as defined in \eqref{eq:w} and $\rho_1 = \max \{ 1+\varepsilon_1 \delta_1, 1+\varepsilon_2 \delta_2 + L_2 \}$. \\
\begin{proof} From \eqref{eq:w} we know 
\begin{equation*}
\begin{aligned}
\|w^{k+1}\|& \leq \| \nabla_u E (u^{k+1},v^{k+1}) +\alpha_1 (p^{k+1}-p^k)\|+\|\nabla_v E (u^{k+1},v^{k+1}) +\alpha_2 (q^{k+1}-q^k)\| \\
& +\| u^k-u^{k+1}\| +\| v^k-v^{k+1}\|
\end{aligned}
\end{equation*}
From the optimality conditons of \eqref{eq:A2.2} and \eqref{eq:A2.4}, we compute
\begin{equation*}
\begin{aligned}
\|w^{k+1}\|& \leq \| \nabla_u E (u^{k+1},v^{k+1}) +\alpha_1 (p^{k+1}-p^k)\|+\|\nabla_v E (u^{k+1},v^{k+1}) +\alpha_2 (q^{k+1}-q^k)\| \\
& +\| u^k-u^{k+1}\| +\| v^k-v^{k+1}\|\\
&= \varepsilon_1 \underbrace{ \|   \nabla G_1 (u^{k+1}) - \nabla G_1 (u^k) \| }_\text{$\leq \delta_1 \|u^{k+1}-u^k \|$}
+ \underbrace{ \|\nabla_u E (u^{k+1},v^{k+1}) - \nabla_u E (u^{k+1},v^k)  \| }_\text{$\leq  L_2 \|v^{k+1}-v^k \|$} \\
&+ \varepsilon_2 \underbrace{ \|  \nabla G_2 (v^{k+1}) - \nabla G_2 (v^k) \| }_\text{$\leq \delta_2 \|v^{k+1}-v^k \|$}
+ \|u^{k+1}-u^k \| + \|v^{k+1}-v^k \| \\
& \leq (1+\varepsilon_1 \delta_1)  \|u^{k+1}-u^k \| +  (1+\varepsilon_2 \delta_2 + L_2) \|v^{k+1}-v^k \| \\
&\leq \rho_1 \|z^{k+1}-z^k\|.
\end{aligned}
\end{equation*}
with $\rho_1 = \max \{ 1+\varepsilon_1 \delta_1, 1+\varepsilon_2 \delta_2 + L_2 \}$. Here we used the  Lipschitz-continuity of $\nabla G_i$ and $\nabla E$.
\end{proof}

Following \cite{BenningBetckeEhrhardtEtAl2017,PALM}, we verify some properties of the limit point set.
Let $\{ z^k \}_{k\in \mathbb{N}}$ and $\{r^k \}_{k\in \mathbb{N}}$ be sequences generated by \eqref{eq:A2}. 
The set of limit points is defined as 
\begin{equation*}
\begin{aligned}
\omega(z^0,r^0) := \Big\{ &(\bar{z},\bar{r}) \in \mathbb{R}^n \times \mathbb{R}^n : \exists \text{ an increasing sequence of integers } \{ k_j \}_{j \in \mathbb{N} }\\
&\text{such that} \lim_{j \to \infty } z^{k_j} = \bar{z}\text{ and } \lim_{j \to \infty } r^{k_j} = \bar{r} \Big\}.
\end{aligned}
\end{equation*}
As in \cite[Definition 5.4, Proposition 5.5]{BenningBetckeEhrhardtEtAl2017}, we are going to assume that $R_i$, $i=1,2$ has locally bounded subgradients.
{\lemma{~\label{lemma_crit}}} Suppose Assumptions \autoref{as:1} and \ref{as:2} hold. Let $\{z^k\}_{k \in \mathbb{N}}$ be a sequence generated by \eqref{eq:A2} which is assumed to be bounded. Let $(\bar{z},\bar{r}) \in \omega(z^0,r^0)$. Then the following assertion holds
\begin{equation}
\lim_{k\to \infty} F(z^{k+1},r^k) = F(\bar{z},\bar{r}) =E(\bar{z}).
\label{eq:crit}
\end{equation}
\begin{proof}
Since $(\bar{z},\bar{r})$ is a limit point of $\{(z^k,r^{k})\}_{k \in \mathbb{N}}$, $\{(z^k,r^{k})\}_{k \in \mathbb{N}}$, there exist subsequences $\{z^{k_j}\}_{j \in \mathbb{N}}$ and $\{r^{k_j}\}_{j \in \mathbb{N}}$ such that $ \lim_{j \to \infty} z^{k_j} = \bar{z}$ and $ \lim_{j \to \infty} r^{k_j} = \bar{r}$, respectively. We immediately obtain
\begin{equation*}
\begin{aligned}
\lim_{j \to \infty} F(z^{k_j},r^{k_j-1}) &= \lim_{j \to \infty} \big\{ E(z^{k_j}) +\alpha_1 D_{R_1}^{p^{k_j-1}}(u^{k_j},u^{k_j-1}) + \alpha_2 D_{R_2}^{q^{k_j-1}}(v^{k_j},v^{k_j-1}) \big\}\\
&= E(\bar{z})
\end{aligned}
\end{equation*}
due to the continuity of $E$ and $\lim_{j \to \infty} D_{R_1}^{p^{k_j-1}}(u^{k_j},u^{{k_j}-1}) =0$ and $\lim_{j \to \infty} D_{R_2}^{q^{k_j-1}}(v^{k_j},v^{k_j-1}) =0$. From the sufficient decrease property we conclude \eqref{eq:crit}.
\end{proof}

{\lemma{(Properties of limit point set).} } The limit point set $w(z^0)$ is a non empty, compact and connected set, the objective function $E$ is constant on $w(z^0)$ and we have $\lim_{k \to \infty} \text{dist} (z^k, w(z^0)) = 0$.
\begin{proof}
This follows steps as in \cite[Lemma 5]{PALM}.
\end{proof}
\newline
\indent To finally prove global convergence of \eqref{eq:A2}, we will use the following Kurdyka-\L{}ojasiewicz property defined  and the result from \cite{PALM}. 
Before recalling the definition, we introduce the notion of distance between any subset $S \subset \mathbb{R}^d$ and any point $x \in \mathbb{R}^d$ defined as
\begin{equation*}
\text{dist}(x,S):= 
\begin{cases} 
\inf \{ \| y - x\| : y\in S\} \qquad &S \neq \emptyset\\
\infty &S=\emptyset 
\end{cases},
\end{equation*} 
where $ \| \cdot \|$ denotes the Euclidean norm.
{\lemma{(Uniformised KL property). \label{KL2}}} Let $\Omega$ be a compact set and let $E: \mathbb{R}^n \times \mathbb{R}^m \to \mathbb{R} \cup \{ \infty \}$ be a proper and l.s.c.\ function. Assume that $E$ is constant on $\Omega$ and satisfy the KL property at each point in $\Omega$. Then there exists $\varepsilon>0$, $\eta>0$ and $\varphi \in C^1((0,\eta))$ that satisfies the same conditions as in Definition KL, such that for all $\bar{u}\in \Omega$ and all $u$ in 
\begin{equation}
\{ u \in \mathbb{R}^n \,\lvert \, \text{dist}(u,\Omega) < \varepsilon \} \cap \{  u \in \mathbb{R}^n \,\lvert \, E(\bar{z}) < E(z)< E(z) + \eta \}
\end{equation}
condition KL is satisfied.
\begin{proof}
Follows from \cite{PALM}.
\end{proof}
\\

With these results we can now show global convergence of \eqref{eq:A2}.

\noindent \begin{proof}[Proof of Theorem 1] By the boundedness assumption on $\{(z^k,r^{k})\}_{k \in \mathbb{N}}$, there exist converging subsequences $\{z^{k_j}\}_{j \in \mathbb{N}}$ and $\{r^{k_j}\}_{j \in \mathbb{N}}$ such that $ \lim_{j \to \infty} z^{k_j} = \bar{z}$ and $ \lim_{j \to \infty} r^{k_j} = \bar{r}$, respectively.
We know from Lemma \autoref{lemma_crit} that \eqref{eq:crit} is satisfied.
\begin{enumerate}
\item 
KL property holds for $E$ and therefore for $E^k$ and we write 
\begin{equation*}
\varphi'\big( F(z^{k},r^{k-1})  - E(\bar{z}) \big) \text{dist}\big(0,\partial F(z^{k},r^{k-1})  \big) \geq 1.
\end{equation*}

From Lemma \autoref{lemmaSubgrBound} we obtain 
\begin{equation*}
\varphi'\big(F(z^{k},r^{k-1})   - E(\bar{z}) \big) \geq \rho_1^{-1} \| z^{k} - z^{k-1} \| ^{-1},
\end{equation*}
and from the concavity of $\varphi$ we know that
\begin{equation*}
\begin{aligned}
\varphi\big(F(z^{k},r^{k-1})   - E(\bar{z}) \big) - \varphi\big((F(z^{k+1},r^{k})   - E(\bar{z})\big)  \\
\geq \varphi' \big(F(z^{k},r^{k-1}) - E(\bar{z}) \big)  \big(F(z^{k},r^{k-1}) -  F(z^{k+1},r^{k}) \big).
\end{aligned}
\end{equation*}
Thus, we obtain
\begin{equation*}
\frac{\varphi \big(F(z^{k},r^{k-1}) - E(\bar{z}) \big) - \varphi\big(F(z^{k},r^{k-1}) - E(\bar{z}) \big)}{F(z^{k},r^{k-1})  -  F(z^{k},r^{k-1})  } 
\geq \rho_{1}^{-1}\| z^{k} - z^{k-1} \| ^{-1}.
\end{equation*}

From \eqref{eq:decrease} with Lemma \autoref{lemmaDecrease} and using the abbreviation 
\begin{equation*}
\varphi^k:=\varphi(F(z^{k},r^{k-1})  - E(\bar{z}) ),
\end{equation*}
it follows 
\begin{equation*}
\begin{aligned}
\frac{ \| z^{k+1} - z^k \|^2} { \| z^{k} - z^{k-1} \| } 
\leq \frac{\rho_1}{\rho_2} (\varphi^k-\varphi^{k+1}).
\end{aligned}
\end{equation*}

Multiplying by $ \| z^{k} - z^{k-1} \| $ and using Young's inequality ($ 2 \sqrt{ab} \leq a+ b $ )

\begin{equation*}
2 \| z^{k+1} - z^k \| \leq \frac{\rho_1} {\rho_2}   (\varphi^k-\varphi^{k+1}) +   \| z^{k} - z^{k-1} \|.
\end{equation*}

Summing up from $k=1,\dots, N$ we get 
\begin{equation*}
\begin{aligned}
\sum_{k=1}^{N} \| z^{k+1} - z^k \| &\leq \frac{\rho_1}{\rho_2}(\varphi^1-\varphi^{N+1})  + \|z^1 - z^0 \| +   \|z^{N+1} - z^N \|\\
 &\leq \frac{\rho_1}{\rho_2}\varphi^1+ \|z^1 - z^0 \|< \infty.
 \end{aligned}
\end{equation*}
In addition we observe that the finite length property implies that the sequence $\{ z^k\}_{k\in \mathbb{N}}$ is a Cauchy sequence and hence is a convergent sequence. For each $z^r$ and $z^s$ with $ s>r>l$ we have

\begin{equation*}
\| z^r - z^s \| = \| \sum_{k=r}^{s-1} z^{k+1} - z^k \| \leq \sum_{k=r}^{s-1} \| z^{k+1} - z^k \|.
\end{equation*} 
\item The proof follows in a similar fashion as in \cite[Lemma 5.9]{BenningBetckeEhrhardtEtAl2017}
\end{enumerate}
\end{proof}

{\remark{(Extension to $d$ blocks).} \label{blocks}}
The analysis described above holds for the general setting of $d$ blocks
\begin{equation}
\min_{\{u_1,\dots,u_d\}} E(u_1, \dots, u_d) + \sum_{i=1}^{n} \alpha_i^k D_{J_i}^{p_i^k}(u_i,u_i^k).
\label{eq:Pblocks}
\end{equation}
The update for each of the d blocks then reads
\begin{equation*}
\begin{aligned}
u_i^{k+1} &= \argmin{u_i} \left\{ E(u_1^{k+1},u_2^{k+1},\dots,u_{i-1}^{k+1},u_i^{k},u_{i+1}^{k}, \dots, u_d^k) + \alpha_i D_{J_i}^{p_i^k}(u_i,u_i^k) \right\} \\
p_i^{k+1} &= p_i^k - \frac{1}{\alpha_i} \Big( \nabla_{u_i} E(u_1^{k+1},u_2^{k+1},\dots,u_{i-1}^{k+1},u_i^{k+1},u_{i+1}^{k}, \dots, u_d^k) \Big).
\end{aligned}
\end{equation*}

\section{Numerical results}
\label{sec:results}
In this section we present numerical results for our joint reconstruction and segmentation model described in \eqref{eq:joint}. We demonstrate its advantages and limitations, as well as a discussion on the parameter choice. In the first part, we focus on bubbly flow segmentation for simulated data. In the second part, we show results for real data acquired at the Cancer Research UK, Cambridge Institute, for tumour segmentation. 
\paragraph{Quality measure.}
To assess the performance of the reconstruction we will compare our solutions $u$ with respect to the groundtruth $u^{gt}$. As quality measure we use the relative reconstuction error (RRE) and  the peak signal to noise ratio (PSNR) defined as
\begin{itemize}
\item $\operatorname{RRE}(u,u^{gt})=\| u^{gt} - u\|_2 / \| u^{gt} \|_2$ 
\item $\operatorname{PSNR}(u,u^{gt})= 10 \log_{10} \left( \frac{\max(u)}{\| u^{gt} - u\|_2 /N}\right)$ 
\end{itemize}
For the segmentation quality, we will use the relative segmentation error (RSE) to compare our segmentations $v$ with respect to the true segmentations $v^{gt}$ 
\begin{itemize}
\item $\operatorname{RSE}(v, v^{gt})= \frac{1}{N} \sum_{i=1}^N \delta_{v^{gt}_i, v_i} $
\end{itemize}
where $N$ is the number of pixels in the image, $\delta$ is the Kronecker delta function that will count the number of mis-classified pixels. 
\newline

\begin{figure}[t!]
\begin{subfigure}{0.23\textwidth}
\centering
 \includegraphics[width=\textwidth]{gt_brain.png}
        \caption{Groundtruth \newline \newline}
        \label{fig:gtbrain1}
    \end{subfigure}
~
\begin{subfigure}{0.23\textwidth}
\centering
 \includegraphics[width=\textwidth]{brain_tvrec.png}
        \caption{TV reconstruction, $\alpha=0.2$,  RRE=0.046, PSNR=24.87 }
        \label{fig:brainTV1}
    \end{subfigure}
    ~
\begin{subfigure}{0.23\textwidth}
\centering
 \includegraphics[width=\textwidth]{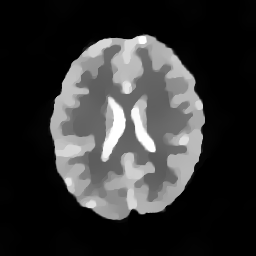}
        \caption{Bregman reconstruction, $\alpha=1$, RRE=0.044, \\PSNR=24.98} 
        \label{fig:brainBreg1}
    \end{subfigure}
~
    \begin{subfigure}{0.23\textwidth}
    \centering
 \includegraphics[width=\textwidth]{brain_jointbregrec1.png}
        \caption{Joint reconstruction,\\$\alpha=0.8$, RRE=0.036,\\ PSNR=26.04} 
        \label{fig:brainJoint1}
    \end{subfigure}

\begin{subfigure}{0.23\textwidth}
\centering
 \includegraphics[width=\textwidth]{map15_col.png}
        \caption{Sampling matrix, 15\% \newline \newline}
        \label{fig:sampling151}
    \end{subfigure}
    ~
            \begin{subfigure}{0.23\textwidth}
        \centering
    \includegraphics[width=\textwidth]{brain_tvseg.png}
        \caption{Segmentation, $\beta=0.001$\\ RSE=0.061  \newline}
        \label{fig:brainTVseg1}
    \end{subfigure}
        ~
\begin{subfigure}{0.23\textwidth}
\centering
 \includegraphics[width=\textwidth]{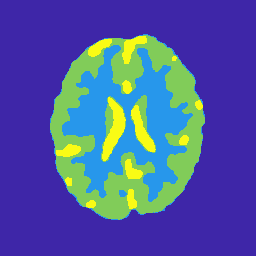}
        \caption{Bregman segmentation, \\$\beta=0.001$, RSE=0.065 \newline}
        \label{fig:brainBregseg1}
    \end{subfigure}
    ~
    \begin{subfigure}{0.23\textwidth}
    \centering
 \includegraphics[width=\textwidth]{brain_jointbregseg1.png}
        \caption{Joint segmentation,\\$\beta=0.001$, $\delta=0.01$\\ RSE=0.057}
        \label{fig:brainJointseg1}
    \end{subfigure} 
    \caption{We consider 15\% of the simulated \textit{k}-space for the brain phantom, where Gaussian noise ($\sigma=0.25$) was added. We compare results for the total variation reconstruction and total variation based Bregman iterative reconstruction and their segmentation in a sequential approach with our joint model. We show that both reconstruction and segmentation are improved.}
    \label{fig:brain_whole}
\end{figure}

Before we present our two applications, we show a more detailed result of the phantom brain in \autoref{fig:brain}. In this example, we show the TV reconstruction \ref{fig:brainTV1}, where the parameter $\alpha$ has been optimised with respect to PSNR and its sequential segmentation \ref{fig:brainTVseg1} with optimal $\beta$ with respect to RSE. In \ref{fig:brainBreg1} and \ref{fig:brainBregseg1} we present Bregman reconstruction and sequential segmentation where the Bregman iteration has been stopped according to the discrepancy principle \autoref{eq:discrepancy} and $\beta$ has been optimised with respect to RSE. These parameter choices for the sequential approaches will be used in the whole paper. \\
In this first result, we clearly see that the joint approach performs much better compared to the separate steps in Figures \ref{fig:brainTV1}, \ref{fig:brainTVseg1} and \ref{fig:brainBreg1}, \ref{fig:brainBregseg1}. Both reconstruction and segmentation are improved and more details are recovered. We refer to \ref{phantoms} for more simulated examples. \textcolor{blue}{}

\subsection{Bubbly flow}
\label{subsec:bubbly}
\begin{figure}[h!]
\centering
\begin{subfigure}{0.3\textwidth}
 \includegraphics[width=\textwidth]{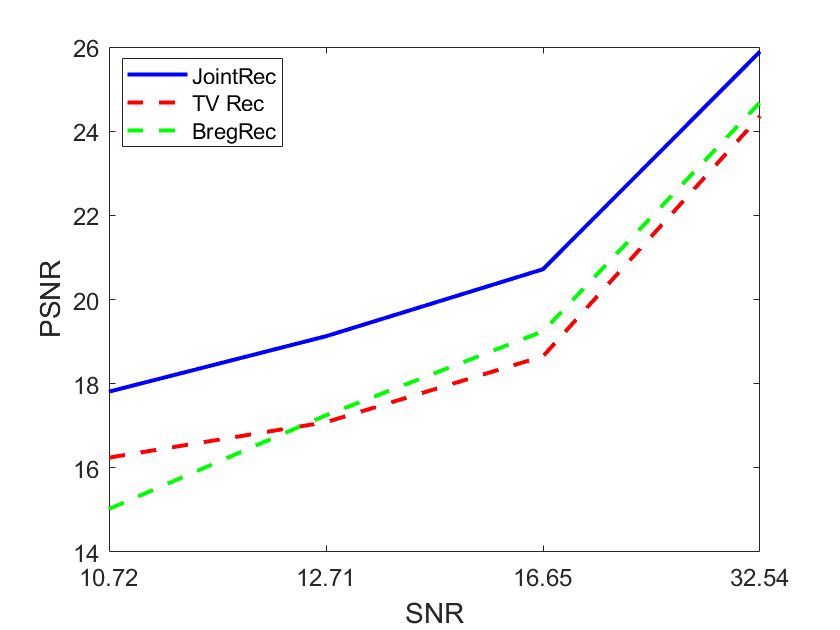}
        \caption{PSNR}
        \label{fig:SNR_PSNR}
    \end{subfigure}
~
    \begin{subfigure}{0.3\textwidth}
 \includegraphics[width=\textwidth]{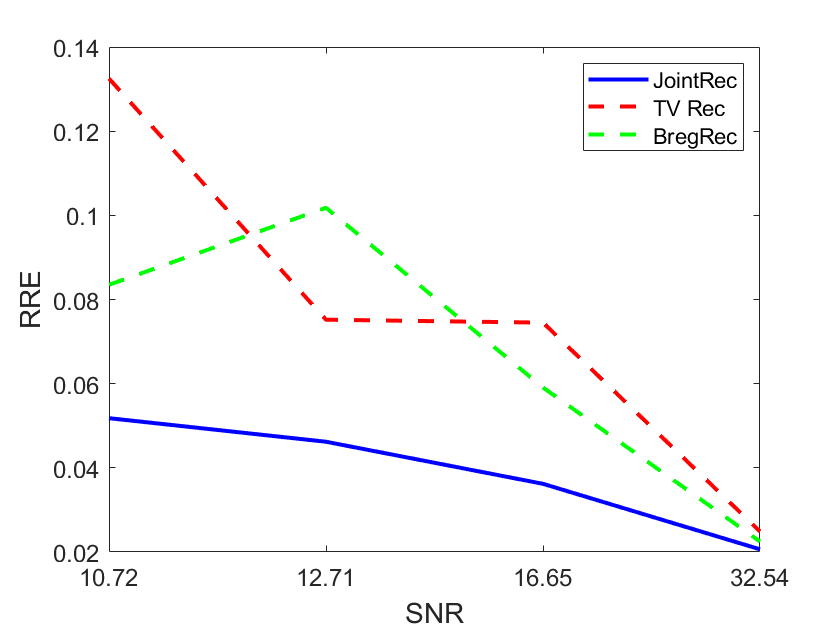}
        \caption{RRE}
        \label{fig:SNR_RRE}
    \end{subfigure}
    ~
    \begin{subfigure}{0.3\textwidth}
 \includegraphics[width=\textwidth]{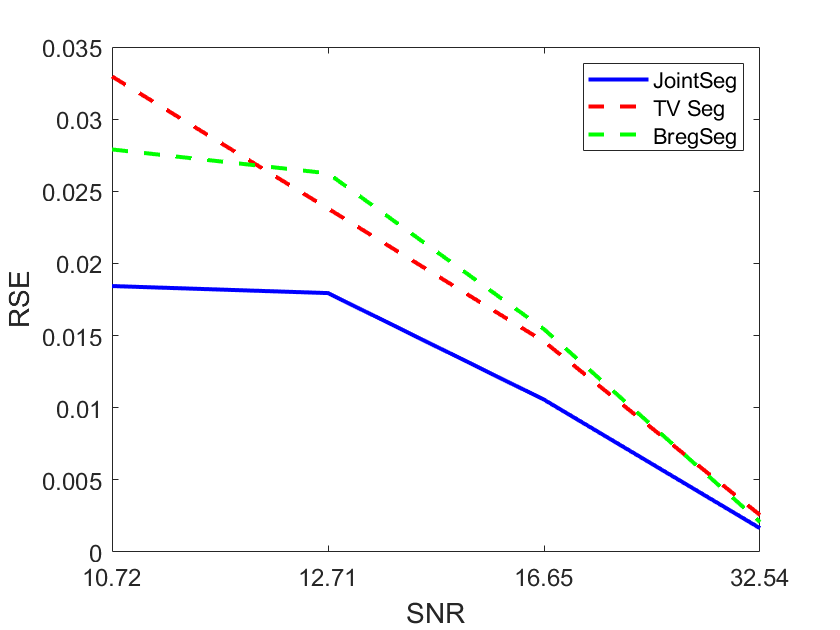}
        \caption{RSE}
        \label{fig:SNR_RSE}
    \end{subfigure}
    \caption{Error plots for different SNR. From left to right, we show the PSNR, RRE and RSE, respectively, for different levels of noise in the measurements. The blue lines represent the error for our joint approach, while the red and green dotted lines are for the sequential TV and sequential Bregman TV approaches. For each SNR, the joint model performs better than the separate methods. This improvement is even more significant for noisier data, which is highly advantageous as in practice we often observe lower SNR. }
    \label{fig:SNR_err}
\end{figure}

\begin{figure}[h!]
\centering
\begin{subfigure}{0.3\textwidth}
 \includegraphics[width=\textwidth]{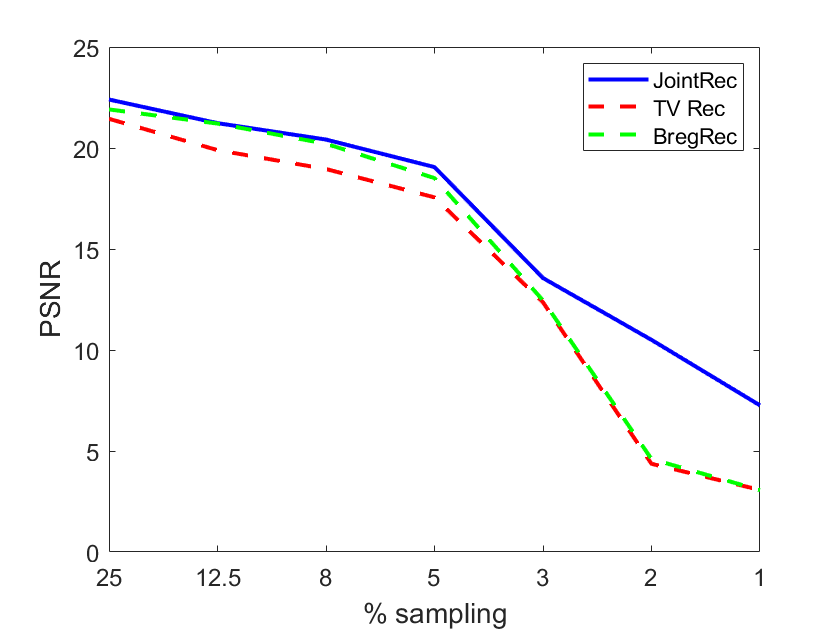}
        \caption{PSNR}
        \label{fig:SAMP_PSNR}
    \end{subfigure}
~
    \begin{subfigure}{0.3\textwidth}
 \includegraphics[width=\textwidth]{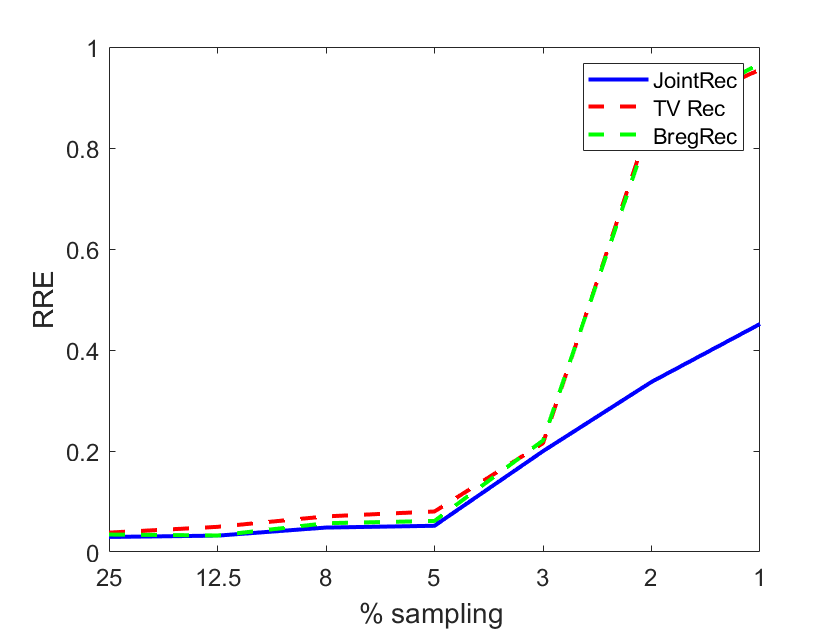}
        \caption{RRE}
        \label{fig:SAMP_RRE}
    \end{subfigure}
~
    \begin{subfigure}{0.3\textwidth}
 \includegraphics[width=\textwidth]{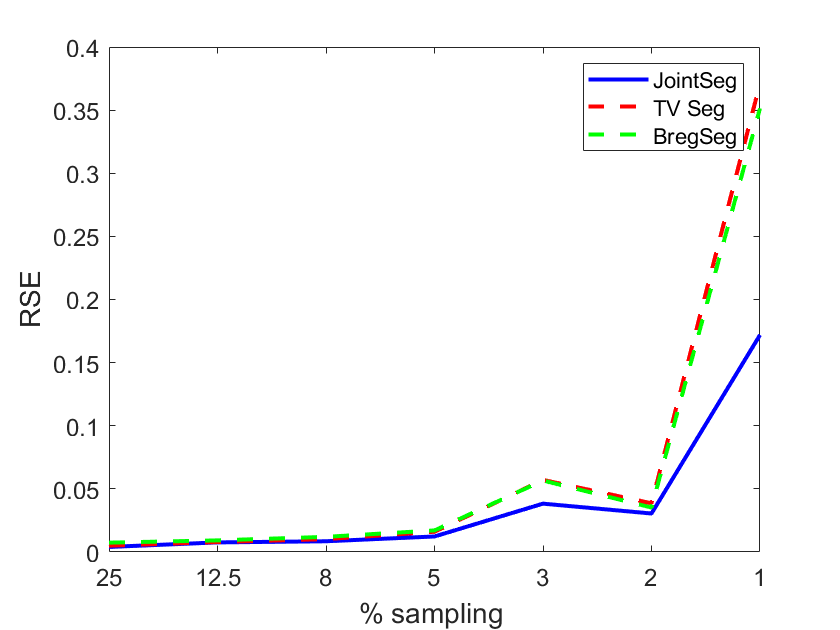}
        \caption{RSE}
        \label{fig:SAMP_RSE}
    \end{subfigure}
    \caption{Error plots varying sampling rate. From left to right, we show the PSNR, RRE and RSE, respectively, for different levels of noise in the measurements. The blue lines represent the error for our joint approach, while the red and green dotted lines are for the sequential TV and sequential Bregman TV approaches. The joint appraoch performs better than the sequential cases. The gain is not very significant for higher sampling rates, but it becomes more important for lower rates, starting from 3\%.}
    \label{fig:SAMP_err}
\end{figure}

The first application considered is the characterisation of bubbly flows using MRI. Bubbly flows are two-phase flow systems of liquid and gas trapped in bubbles, which are common in industrial applications such as bioreactors \cite{Chalmers1994} and hydrocarbon processing units \cite{Deckwer1992}. MRI has been successfully used to characterise the bubble size distribution \cite{Holland2012,Tayler2012} and the liquid velocity field of bubbly flows \cite{HollandMalioutovBlakeEtAl2010,TaylerHollandSedermanEtAl2012}; these properties govern the heat and mass transfer between the bubbles and the liquid which ultimately determine the efficiency of these industrial systems. However, when studying fast flowing systems, the acquisition time for fully sample \textit{k}-space is too long to resolve the temporal changes; the most common method of breaking the temporal resolution barrier is through under-sampling. It is therefore critical to develop reconstruction techniques for highly under-sampled \textit{k}-space data for the accurate reconstruction of the MRI images which would be subsequently used in calculating the bubble size distribution or in studying the hydrodynamics of the system.\\
\begin{figure}[t]
\centering
\begin{subfigure}{0.23\textwidth}
\centering
 \includegraphics[width=\textwidth]{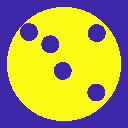}
        \caption{Groundtruth \newline  \newline}
        \label{fig:gt}
    \end{subfigure}
~
\begin{subfigure}{0.23\textwidth}
\centering
 \includegraphics[width=\textwidth]{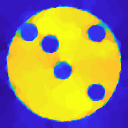}
        \caption{$\operatorname{TV}$ reconstruction, $\alpha=0.1$, RRE=0.081, PSNR=18.42 }
        \label{fig:bubblyTV}
    \end{subfigure}
    ~
\begin{subfigure}{0.23\textwidth}
\centering
 \includegraphics[width=\textwidth]{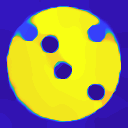}
        \caption{Bregman reconstruction, $\alpha=2$, RRE=0.069, PSNR=18.83} 
        \label{fig:bubblyBreg}
    \end{subfigure}
~
    \begin{subfigure}{0.23\textwidth}
    \centering
 \includegraphics[width=\textwidth]{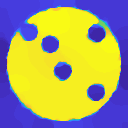}
        \caption{Joint reconstruction, \\$\alpha=0.8$,RRE=0.058, PSNR=20.7105 } 
        \label{fig:bubblyJoint}
    \end{subfigure}

\begin{subfigure}{0.23\textwidth}
\centering
 \includegraphics[width=\textwidth]{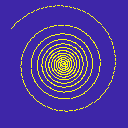}
        \caption{Sampling matrix, 8\% \newline }
        \label{fig:sampling}
    \end{subfigure}
    ~
            \begin{subfigure}{0.23\textwidth}
        \centering
    \includegraphics[width=\textwidth]{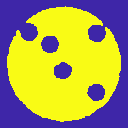}
        \caption{Segmentation, \\$\beta=0.001$, RSE=0.0093}
        \label{fig:bubblyTVseg}
    \end{subfigure}
        ~
\begin{subfigure}{0.23\textwidth}
\centering
 \includegraphics[width=\textwidth]{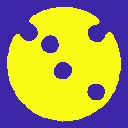}
        \caption{Bregman segmentation, $\beta=0.001$, RSE=0.017}
        \label{fig:bubbleBregseg}
    \end{subfigure}
    ~
    \begin{subfigure}{0.23\textwidth}
    \centering
 \includegraphics[width=\textwidth]{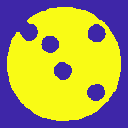}
        \caption{Joint segmentation, $\beta=0.001$, $\delta=1$, RSE=0.0102} 
        \label{fig:bubblyJointseg}
    \end{subfigure} 
    \caption{Results of  the $\operatorname{TV}$ reconstruction and  Bregman iterative reconstruction and their segmentation in the sequential approach are compared with our joint model. Both MSE and SSIM are improved in the joint approach. The data was corrupted with Gaussian noise with $\sigma=0.35$.}
    \label{fig:bubblyflow}
    \end{figure}
\indent We apply our joint reconstruction and segmentation approach to simulated bubbly flow imaging.
In \autoref{fig:bubblyflow} we present some results for synthetic data, where \autoref{fig:gt} represents the groundtruth magnitude image, from which we simulate its \textit{k}-space following the forward model described in \eqref{eq:generalforward}. From the full \textit{k}-space we collect 8\% of the samples using the sampling matrix in \autoref{fig:sampling} and we corrupt the data with Gaussian noise of standard deviation $\sigma = 0.35$. In \autoref{fig:bubblyTV} and \ref{fig:bubblyTVseg} we show the results for the total variation regularised reconstruction and its segmentation performed sequentially. In the same sequential way, we show the results for the Bregman iterative regularization in \autoref{fig:bubblyBreg} and \ref{fig:bubbleBregseg}. In the last column in \autoref{fig:bubblyJoint} and \ref{fig:bubblyJointseg}, we finally show the results for our joint approach. Although the TV and the Bregman approaches are already quite good, we can see that both RRE and PSNR are improved using our model in the reconstruction and the segmentation. Smaller details, such as the top right bubble contour, are better detected when solving the joint problem. As the goal of the bubbly flow application is to detect bubble size distribution, this improvement is really advantageous. \\ 
\newline
\indent We tested the robustness of our approach by corrupting the data with different signal to noise ratio (SNR) and by considering different amount of sampling. In \autoref{fig:snr} we show in the top row the reconstructions obtained with the joint model for different SNR (which corresponds to different standard deviation $\sigma$) and in the bottom row the corresponding segmentation obtained by the joint approach. To complement this information, we show in \autoref{fig:SNR_err} how the PSNR, RRE and RSE are affected, for the joint approach (blue lines) and for the separate approaches, TV (red dotted lines) and Bregman TV (green dotted lines). As expected, with the SNR increasing the error decreases. We can see that the joint approach performs better than the sequential approach for any SNR. The improvement is even more significant for very noisy data. As in practice we often observe high levels of noise, the joint approach is able to takle this problem better than the traditional sequential approaches.\\ 
\indent It is also interesting to investigate how the joint approach performs with very low undersampling rates. In Figure \ref{fig:sampling} we show joint reconstructions (top row) and corresponding segmentations (bottom row) for decreasing sampling rates. We can see that up to 5\% results are still very good. Using 3 and 2\% of the samples the results are less clean but it is possible to identify the main structures. In contrast, 1\% sampling is not enough to retrieve a good image reconstruction and consequently its segmentation. In \autoref{fig:SAMP_err}, we plot PSNR , RRE and RSE  for different sampling rates. The blue lines represent the error for our joint approach, while the red and green dotted lines are for the sequential TV and sequential Bregman TV approaches. We can see that up to 5\% sampling the error measures do not change significantly. However, for lower rates, the improvement is more significant. This is highly beneficial for the bubbly flow application as increasing the temporal resolution is really important to keep track of the gas flowing in the pipe.

\begin{figure}[t!]
\centering
\begin{subfigure}{0.23\textwidth}
\centering
 \includegraphics[width=\textwidth]{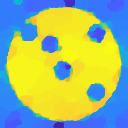}
        \label{fig:}
    \end{subfigure}
~
\begin{subfigure}{0.23\textwidth}
\centering
 \includegraphics[width=\textwidth]{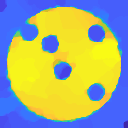}
        \label{fig:}
    \end{subfigure}
    ~
\begin{subfigure}{0.23\textwidth}
\centering
 \includegraphics[width=\textwidth]{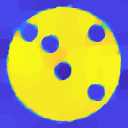}
        \label{fig:}
    \end{subfigure}
~
    \begin{subfigure}{0.23\textwidth}
    \centering
 \includegraphics[width=\textwidth]{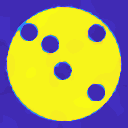}
        \label{fig:}
    \end{subfigure}

\begin{subfigure}{0.23\textwidth}
\centering
 \includegraphics[width=\textwidth]{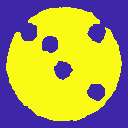}
        \caption{SNR=10.56, $\sigma$=0.70}
        \label{fig:}
    \end{subfigure}
    ~
            \begin{subfigure}{0.23\textwidth}
        \centering
    \includegraphics[width=\textwidth]{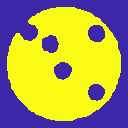}
        \caption{SNR=12.69, $\sigma$=0.56}
        \label{fig:}
    \end{subfigure}
        ~
\begin{subfigure}{0.23\textwidth}
\centering
 \includegraphics[width=\textwidth]{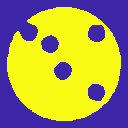}
        \caption{SNR=16.68, $\sigma$=0.35}
        \label{fig:}
    \end{subfigure}
    ~
    \begin{subfigure}{0.23\textwidth}
    \centering
 \includegraphics[width=\textwidth]{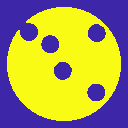}
        \caption{SNR=32.83, $\sigma$=0.06}
        \label{fig:}
    \end{subfigure} 
\caption{Top row: reconstructions obtained by the joint model with different SNR. Bottom row: corresponding segmentations.}
\label{fig:snr}
\end{figure}

\begin{figure}[t]
\centering
 \includegraphics[width=\textwidth]{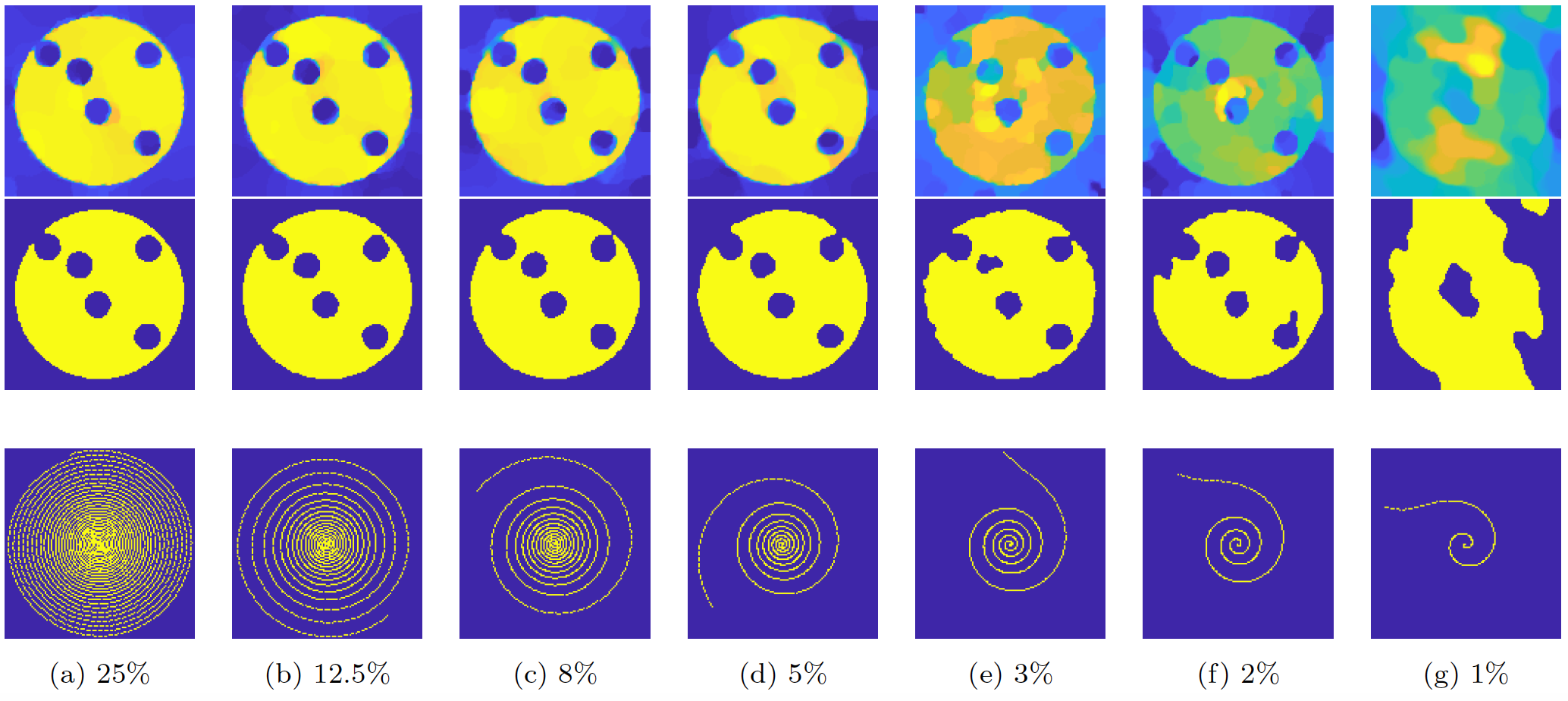}

\caption{Top row: reconstructions obtained by the joint model with different sampling rates. Bottom row: corresponding segmentations. The joint reconstruction and segmentation is able to detect the main structures down  to 5\% of the samples. Up to 2\% the results are less clean but still acceptable. Using only 1\% of the data is not enough to produce the image and segmentation. }
\label{fig:snr}
\end{figure}

\begin{figure}[t]
\centering
\begin{subfigure}{0.23\textwidth}
 \includegraphics[width=\textwidth]{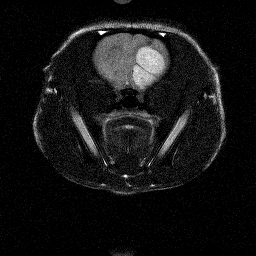}
        \caption{Zero-filled reconstruction}
        \label{fig:tum_zero}
    \end{subfigure}
~
\begin{subfigure}{0.23\textwidth}
 \includegraphics[width=\textwidth]{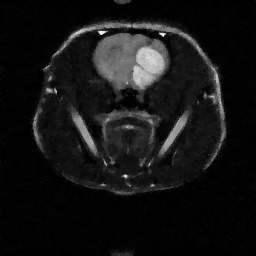}
        \caption{TV reconstruction \\$\alpha=0.01$}
        \label{fig:tum_tvrec}
    \end{subfigure}
    ~
\begin{subfigure}{0.23\textwidth}
 \includegraphics[width=\textwidth]{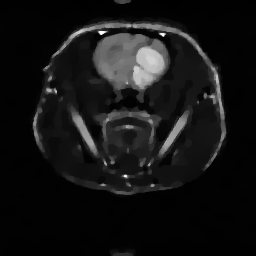}
        \caption{Bregman reconstruction \\$\alpha=1$}
        \label{fig:tum_bregrec}
    \end{subfigure}
~
    \begin{subfigure}{0.23\textwidth}
 \includegraphics[width=\textwidth]{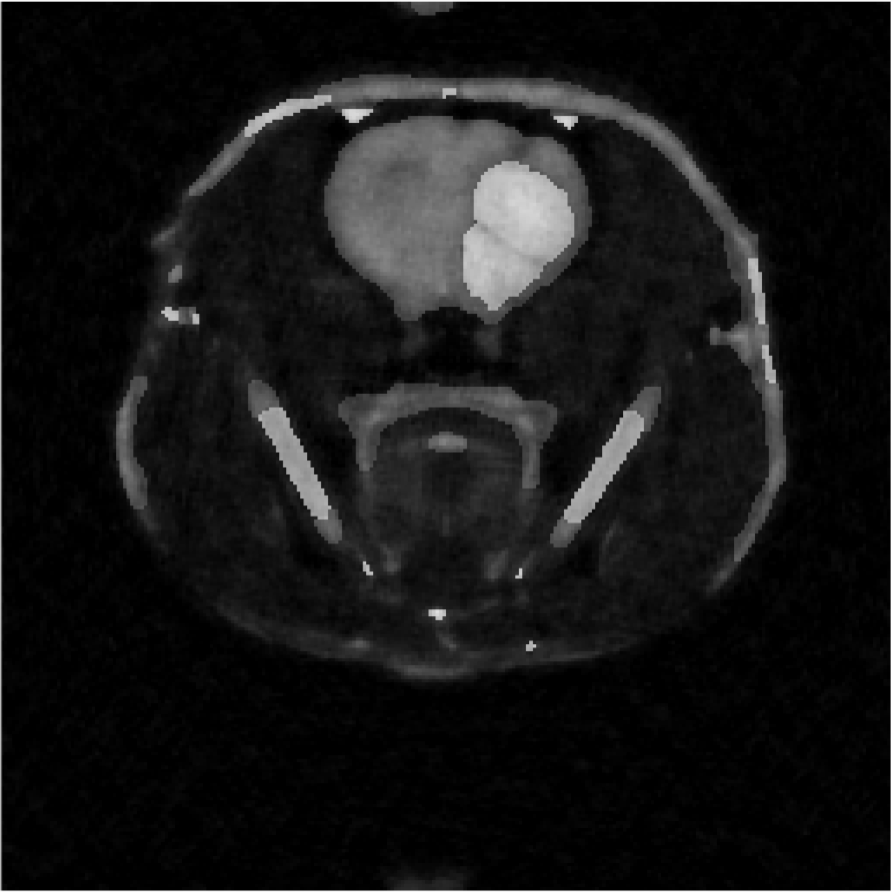}
        \caption{Joint reconstruction \\ $\alpha=0.5$}
        \label{fig:tum_jointrec}
    \end{subfigure}
    
           \begin{subfigure}{0.23\textwidth}
    \includegraphics[width=\textwidth]{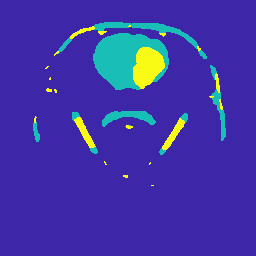}
        \caption{Segmentation \newline }
        \label{fig:tum_zero_seg}
    \end{subfigure}
~
        \begin{subfigure}{0.23\textwidth}
    \includegraphics[width=\textwidth, height=3.62cm]{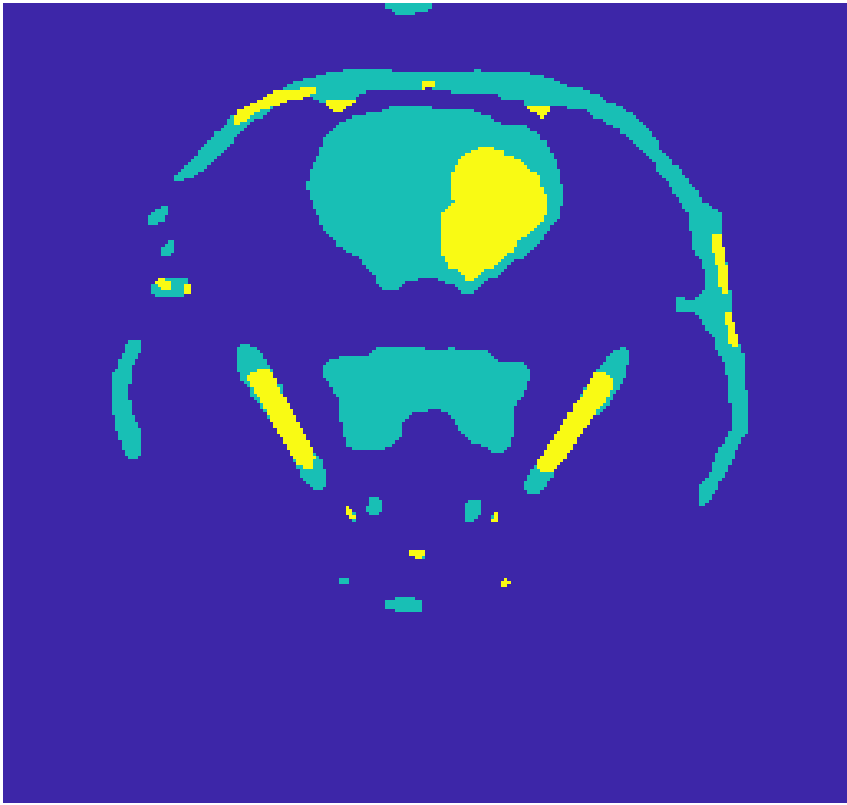}
        \caption{Segmentation \\ $\beta=0.07$ }
        \label{fig:tum_tvseg}
    \end{subfigure}
~
\begin{subfigure}{0.23\textwidth}
 \includegraphics[width=\textwidth]{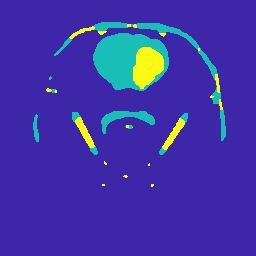}
        \caption{Bregman segmentation \\$\beta=0.07$}
        \label{fig:tum_bregseg}
    \end{subfigure}
    \begin{subfigure}{0.23\textwidth}
    ~
 \includegraphics[width=\textwidth, height=3.62cm]{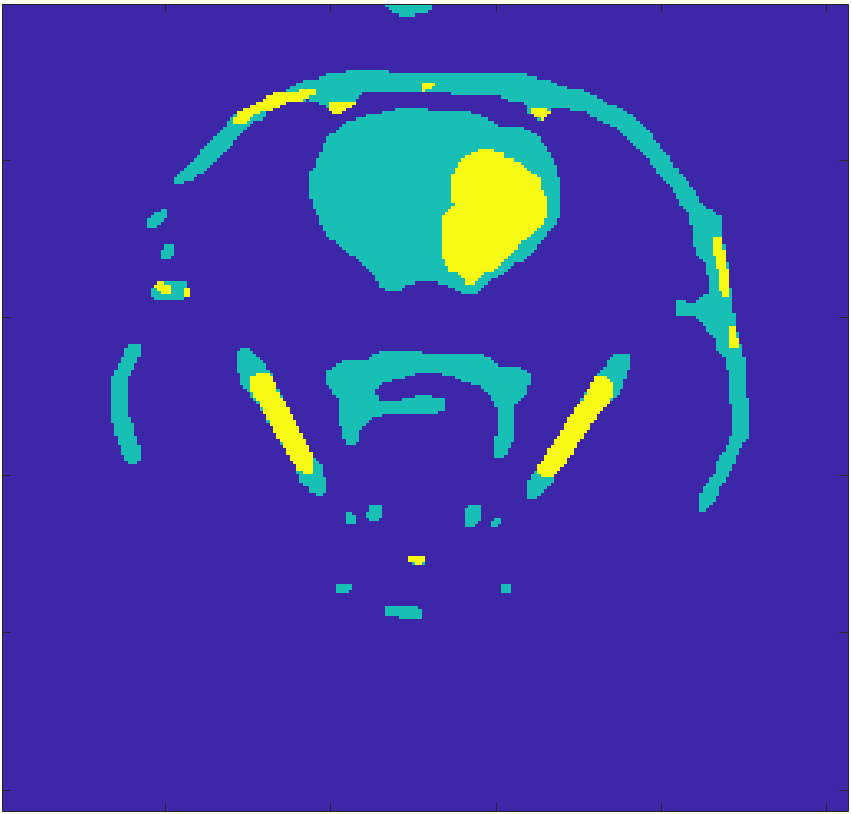}
        \caption{Joint segmentation \\$\beta=0.01$, $\delta=0.01$ }
        \label{fig:tum_jointseg}
    \end{subfigure} 
    \caption{Reconstructions and segmentation for real MRI data. We select 15\% of the samples using a spiral mask. The image show a rat brain bearing a tumour (brighter region). The zero-filled reconstruction \ref{fig:tum_zero} and the TV regularised reconstruction \ref{fig:tum_tvrec} are shown together with their sequential segmentation \ref{fig:tum_zero_seg} and \ref{fig:tum_tvseg} respectively. In the last column \ref{fig:tum_jointrec} and \ref{fig:tum_jointseg} we show the results for our model. The parameter $\alpha$ for the TV reconstruction and for the joint reconstruction has been chosen such that it achieves visually optimal in the sense that it resolve all the details (e.g. the darker line cutting the tumour transversally).}
    \label{fig:tum}
\end{figure}
\begin{figure}[h!]
\centering
\begin{subfigure}{0.2\textwidth}
 \includegraphics[width=\textwidth]{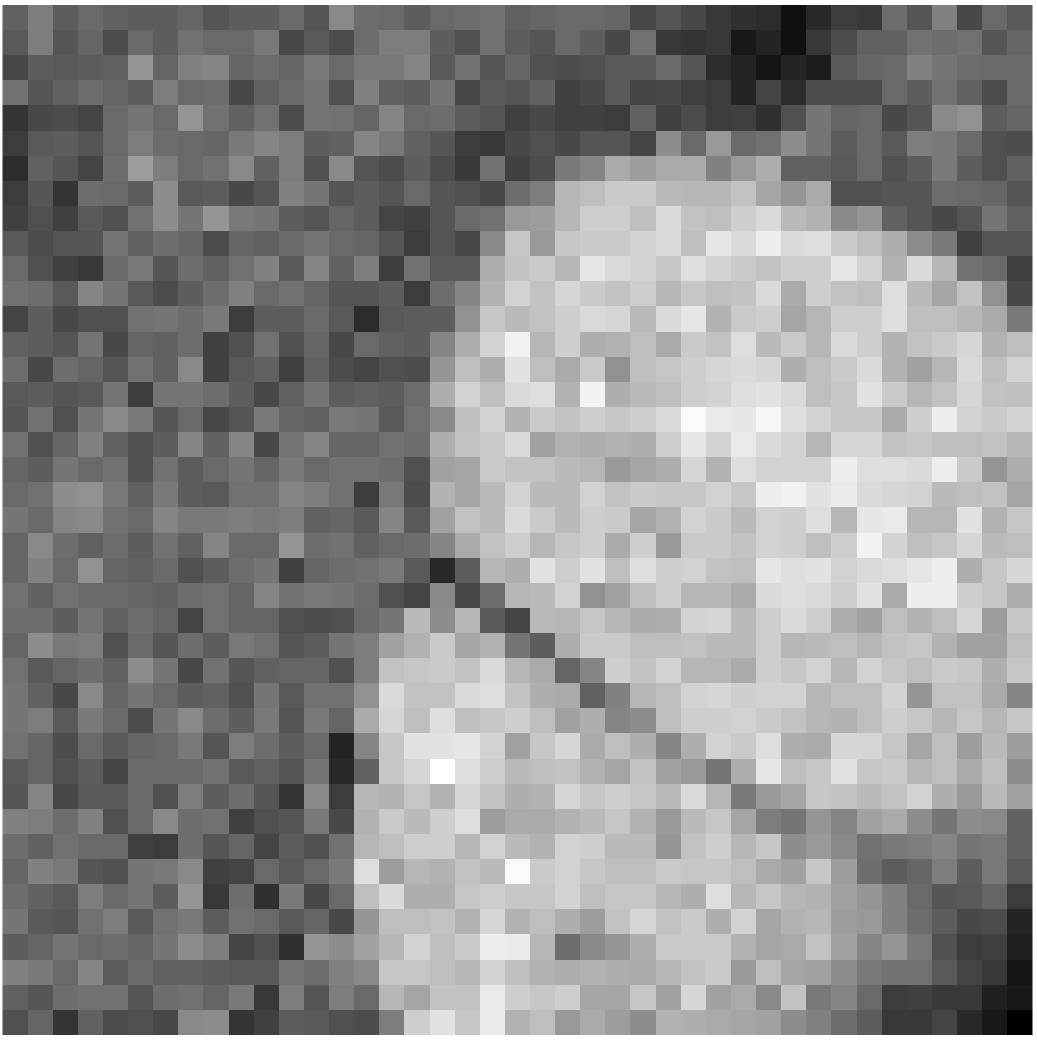}
        \caption{Zero-filled reconstruction }
        \label{fig:}
    \end{subfigure}
~
\begin{subfigure}{0.2\textwidth}
 \includegraphics[width=\textwidth]{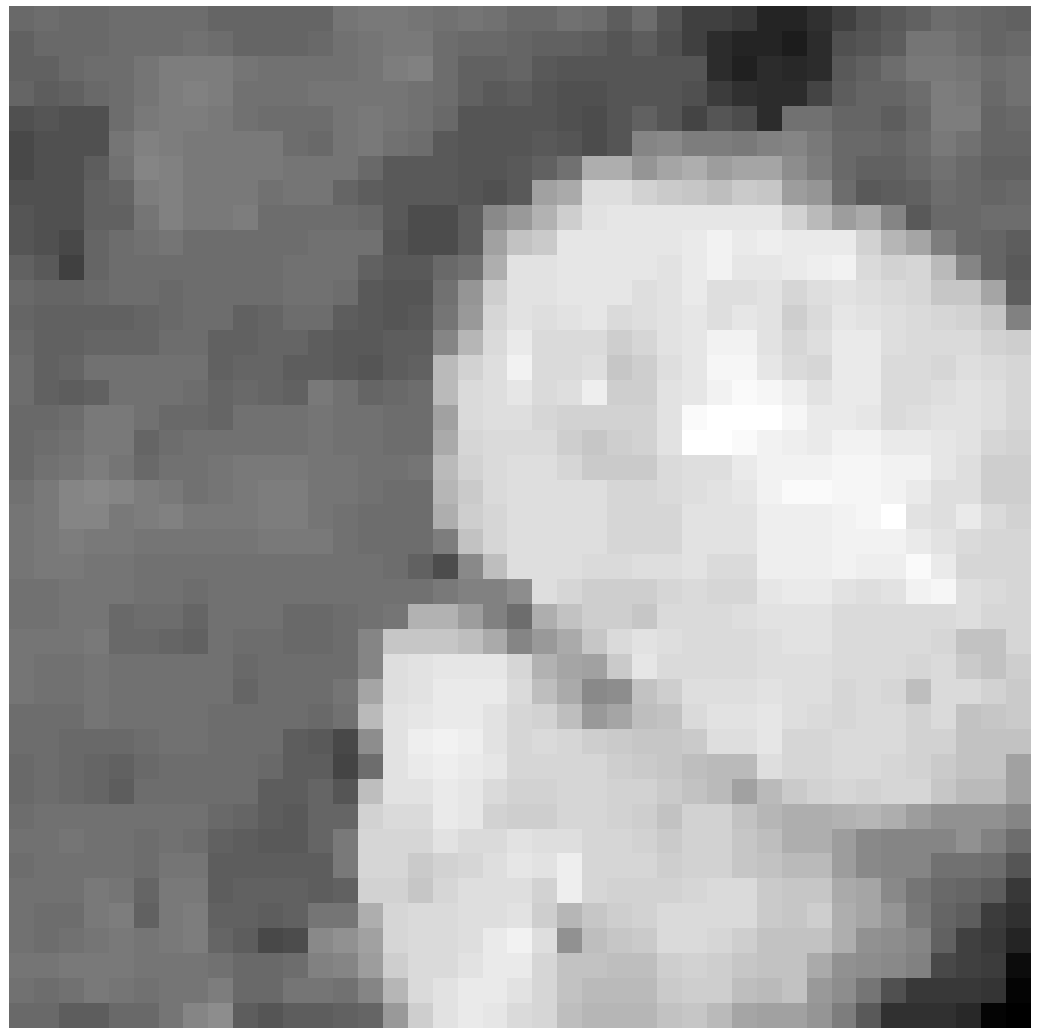}
        \caption{TV reconstruction \newline}
        \label{fig:}
    \end{subfigure}
    ~
    \begin{subfigure}{0.2\textwidth}
 \includegraphics[width=\textwidth]{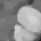}
        \caption{Bregman reconstruction}
        \label{fig:}
    \end{subfigure}
~
    \begin{subfigure}{0.2\textwidth}
 \includegraphics[width=\textwidth]{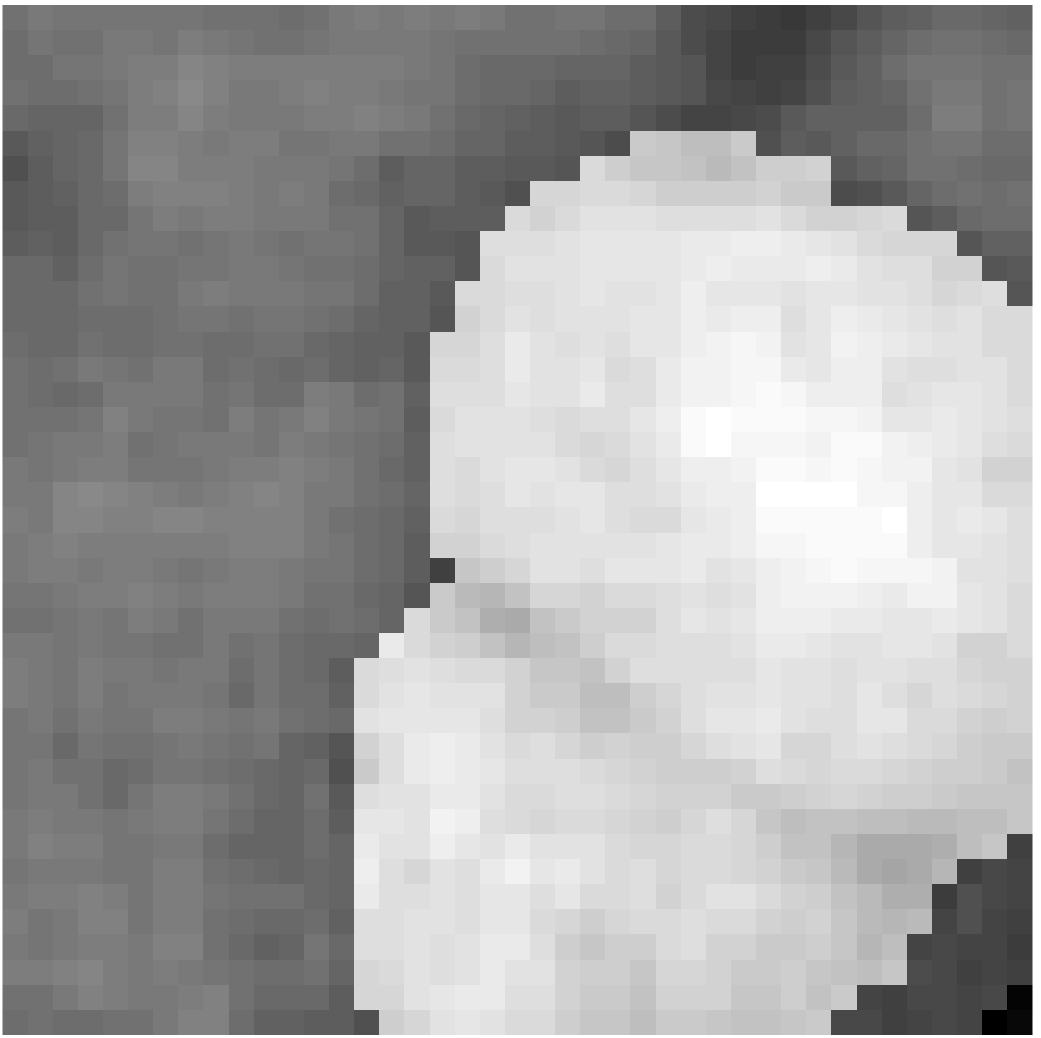}
        \caption{Joint reconstruction \newline}
        \label{fig:}
    \end{subfigure}
    
        \begin{subfigure}{0.2\textwidth}
    \includegraphics[width=\textwidth]{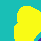}
        \caption{Segmentation  \newline}
        \label{fig:}
    \end{subfigure}
~
       \begin{subfigure}{0.2\textwidth}
    \includegraphics[width=\textwidth ,height=3.15cm]{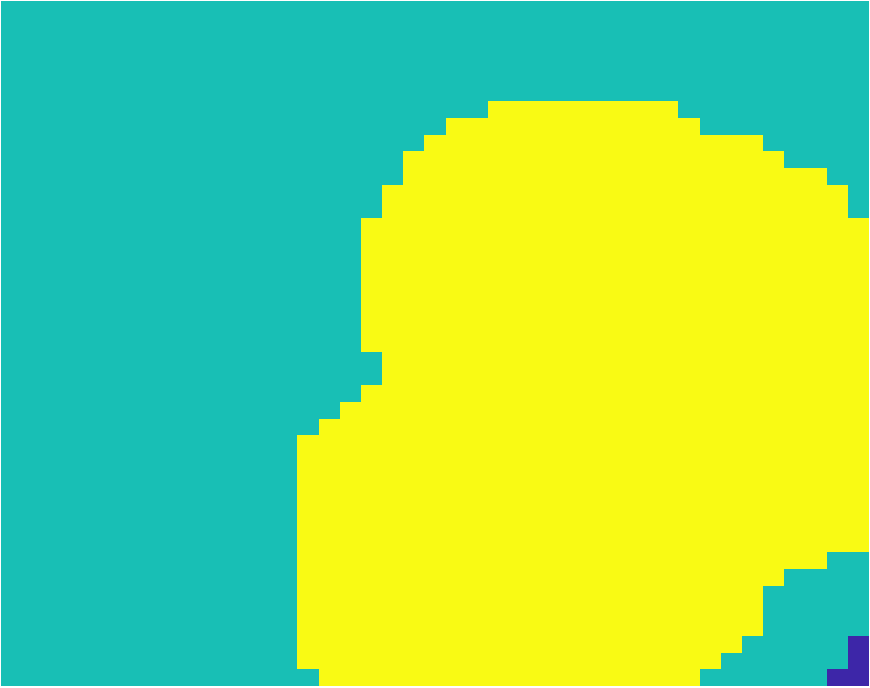}
        \caption{Segmentation \newline}
        \label{fig:}
    \end{subfigure}
    ~
       \begin{subfigure}{0.2\textwidth}
    \includegraphics[width=\textwidth , height=3.15cm]{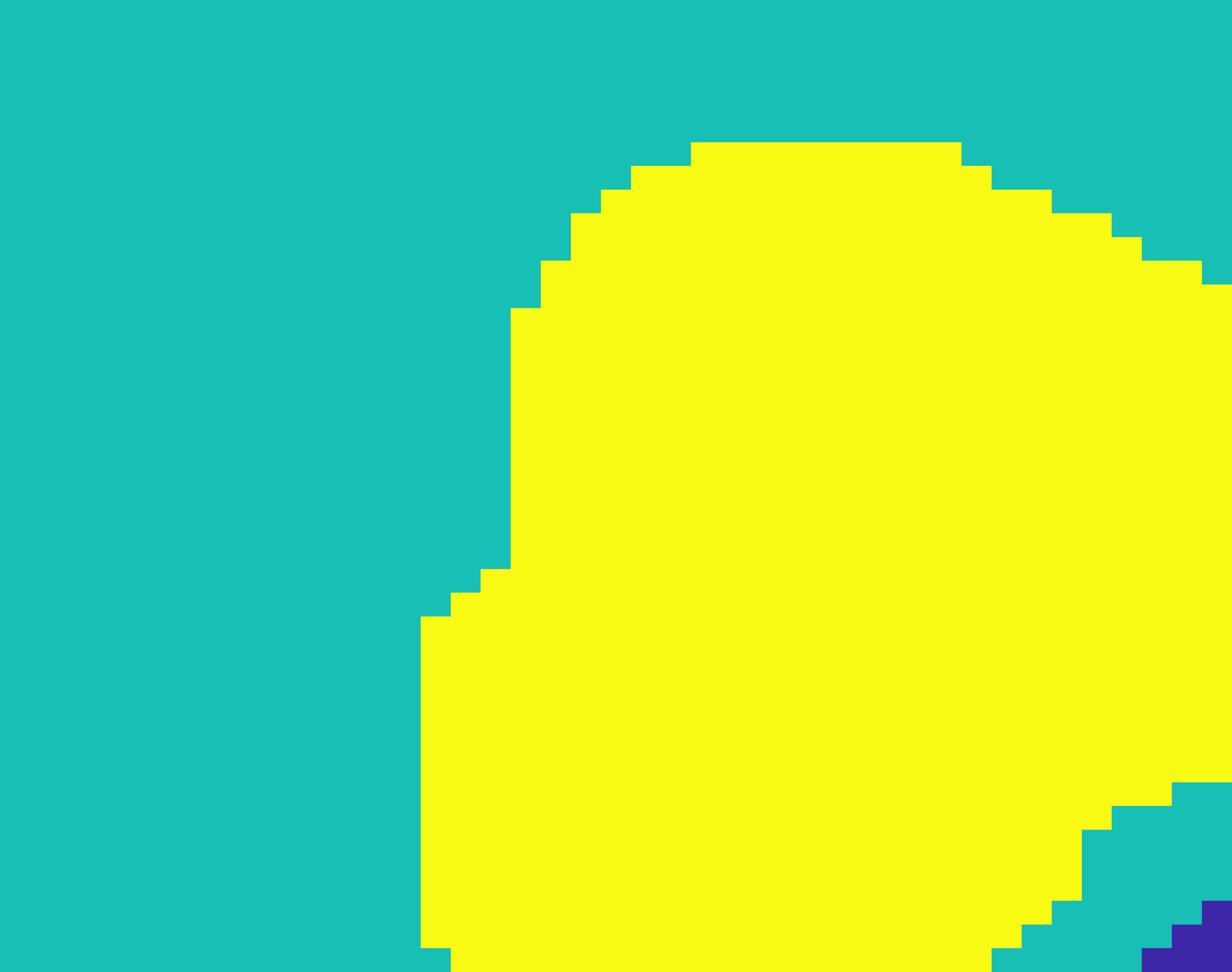}
        \caption{Bregman segmentation}
        \label{fig:}
    \end{subfigure}
~
    \begin{subfigure}{0.2\textwidth}
 \includegraphics[width=\textwidth ,height=3.15cm]{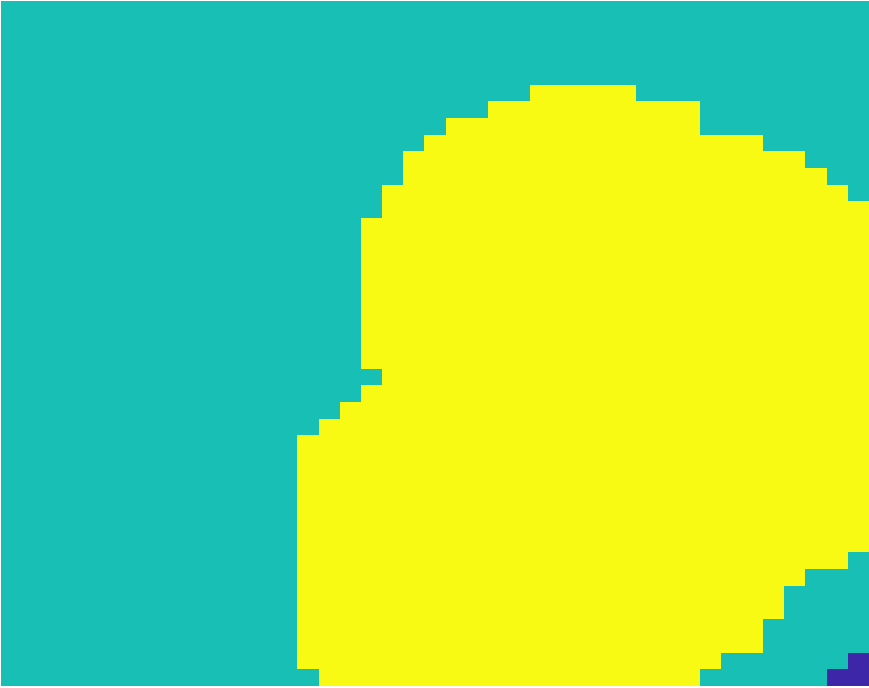}
        \caption{Joint segmentation  \newline}
        \label{fig:}
    \end{subfigure} 
    \caption{Zoomed section on the tumor for the different approaches.}
    \label{fig:tum_zoom}
\end{figure}
\subsection{Cancer imaging}

In this subsection, we illustrate the performance of the joint model for real cancer data. At the Cancer Research UK, Cambridge Institute, researchers acquire every day a huge amount of MRI scans to assess tumour progression and response to therapy \cite{Rodrigues2014}. For this reason, it is very convenient to have fast sampling through compressed sensing, and automatic segmentation methods. Furthermore, reconstructions with enhanced edges are advantageous to facilitate clinical analysis. \\
\indent Here we show our results for MRI data of a rat bearing a glioblastoma. The MR image represents the rat head where the brain is the gray area in the top half of the image. Inside this gray region, a tumour is clearly visible appearing as a brighter area. For this experiment, we acquired the full \textit{k}-space and present the zero-filling reconstruction in \autoref{fig:tum_zero} and the sequential segmentation in \autoref{fig:tum_zero_seg}. As discussed already in the previous section, the zero-filled reconstruction presents noise and artefact which may complicate the segmentation. We want to show that the compressed sensing approach and in particular the joint model can improve this reconstruction. Given the full \textit{k}-space, we select 15\% of the samples using a spiral mask. In \autoref{fig:tum_tvrec}, \ref{fig:tum_tvseg} and \autoref{fig:tum_bregrec},  \ref{fig:tum_bregseg} we show the results for the sequential approaches. In \autoref{fig:tum_jointrec} and \ref{fig:tum_jointseg} we show the joint reconstruction and the joint segmentation obtained for the same data. The regularised approaches perform better that the zero-filled reconstruction, producing less noisy results. However, our joint model is able to produce a cleaner reconstruction where the edges that defines the tumour and the brain are very well detected. In \autoref{fig:tum_zoom}, we show a zoomed section where it is easy to assess that the joint model tackle the noise and detect the region of interest. We can see that we are able to improve the reconstruction and automatically identify the tumour in the brain.  The degree of enhancement of the edges in the reconstruction is controllable by the parameter $\delta$ in the model \eqref{eq:joint}. In the next subsection we present a discussion on how to tune this parameter. 

\subsection{Parameter choice rule}

\begin{figure}[t!]
\centering
 \begin{subfigure}{0.2\textwidth}
 \includegraphics[width=\textwidth]{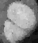}
        \caption{$\alpha = 0.001 $}
        \label{fig:alpha0001}
    \end{subfigure}
~
\begin{subfigure}{0.2\textwidth}
 \includegraphics[width=\textwidth]{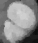}
        \caption{$\alpha = 0.01 $ }
        \label{fig:alpha001}
    \end{subfigure}
    ~
\begin{subfigure}{0.2\textwidth}
 \includegraphics[width=\textwidth]{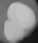}
        \caption{$\alpha = 0.1 $}
        \label{fig:alpha01}
    \end{subfigure}
 \caption{The parameter $\alpha$ balances the data fidelity term and the total variation regularisation for the reconstruction. Smaller values of $\alpha$  produce a reconstruction closer to the data fitting term, hence less smooth as in \ref{fig:alpha0001}. As $\alpha$ increases in \ref{fig:alpha001} the solution gets smoother and less noisy.  Finally for large values it tends to become more piecewise constant as in \ref{fig:alpha01}.} 
 \label{fig:alphas}
\end{figure}

\begin{figure}[t!]
\centering
\begin{subfigure}{0.2\textwidth}
 \includegraphics[width=\textwidth]{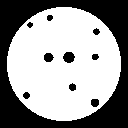}
        \caption{$\beta = 0.1 $}
        \label{fig:beta01}
    \end{subfigure}
~
\begin{subfigure}{0.2\textwidth}
 \includegraphics[width=\textwidth]{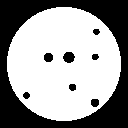}
        \caption{$\beta = 1 $ }
        \label{fig:beta1}
    \end{subfigure}
~
    \begin{subfigure}{0.2\textwidth}
 \includegraphics[width=\textwidth]{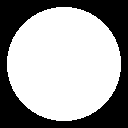}
        \caption{$\beta = 3 $}
        \label{fig:beta3}
    \end{subfigure}
 \caption{The parameter $\beta$ determines the scale of the objects that we are segmenting. Smaller values of $\beta$ can detect smaller objects \ref{fig:beta01}, which are lost for intermediate values \ref{fig:beta1}. Finally very large values only detect main structures  \ref{fig:beta3}. }
 \label{fig:betas}
\end{figure}

\begin{figure}[t!]
\centering
\begin{subfigure}{0.2\textwidth}
 \includegraphics[width=\textwidth]{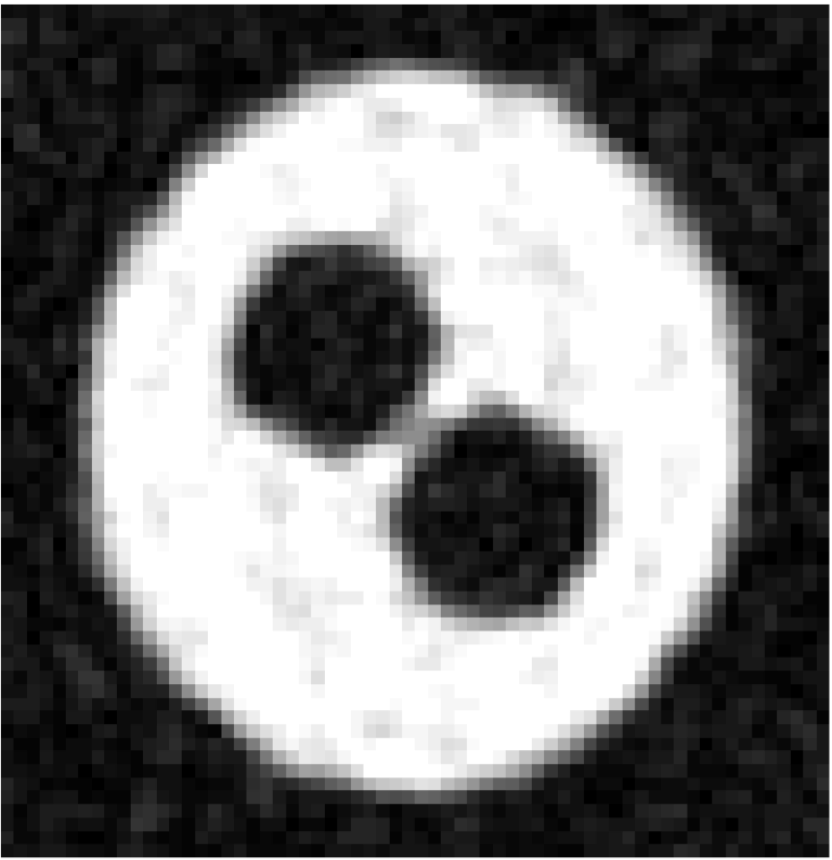}
        \caption{$\delta = 0 $}
        \label{fig:}
    \end{subfigure}
~
\begin{subfigure}{0.2\textwidth}
 \includegraphics[width=\textwidth]{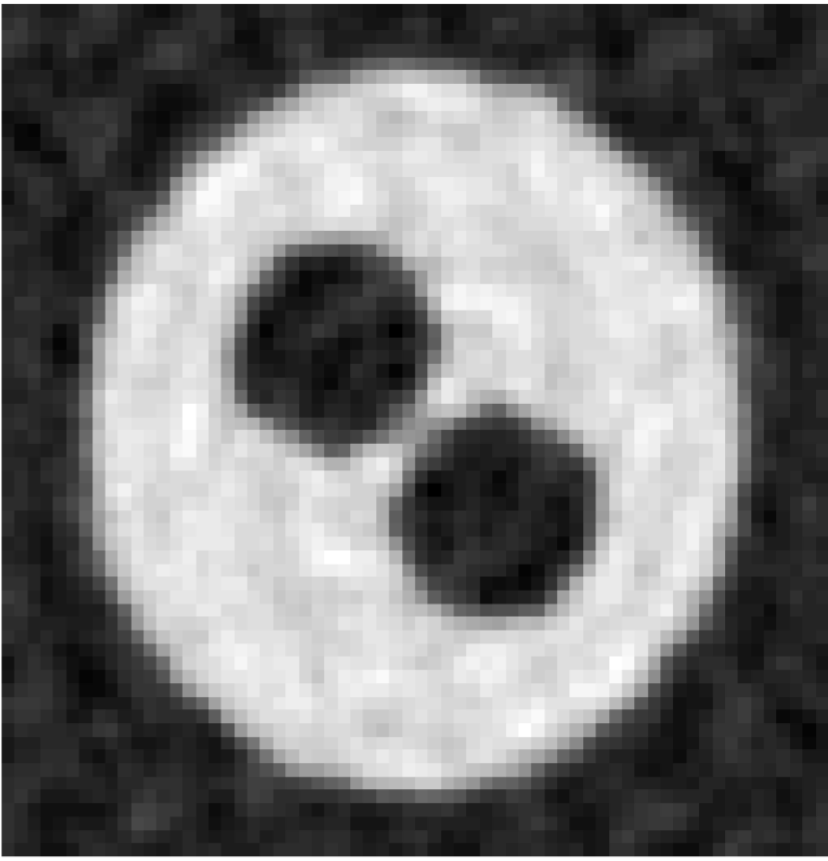}
        \caption{$\delta = 0.1 $ }
        \label{fig:}
    \end{subfigure}
~
    \begin{subfigure}{0.2\textwidth}
 \includegraphics[width=\textwidth]{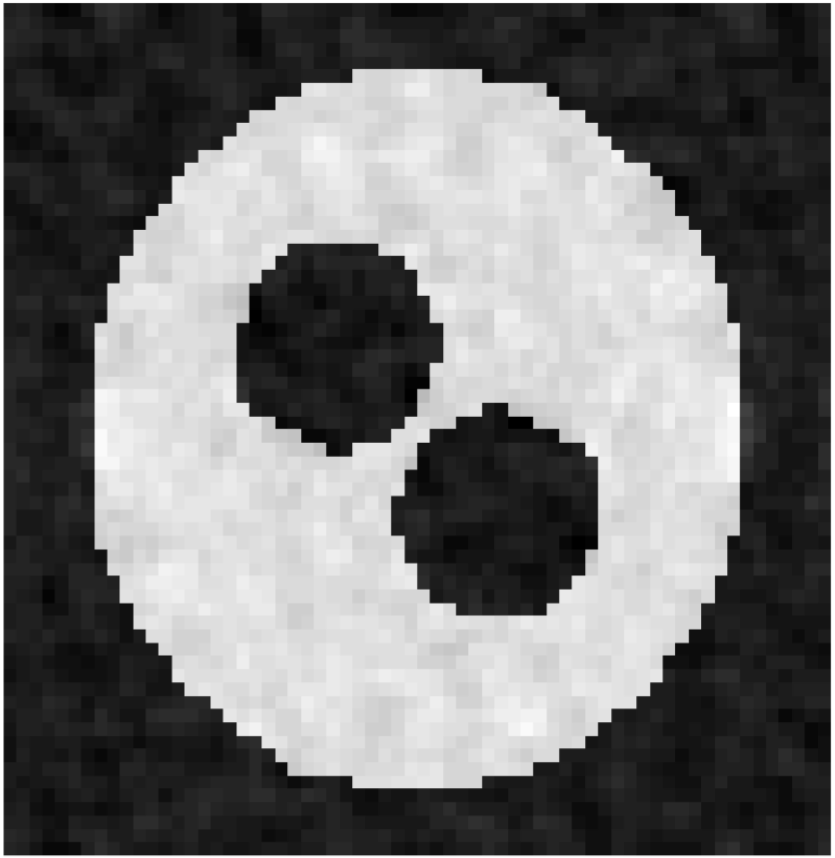}
        \caption{$\delta = 1 $}
        \label{fig:}
    \end{subfigure}
    ~
    \begin{subfigure}{0.2\textwidth}
 \includegraphics[width=\textwidth]{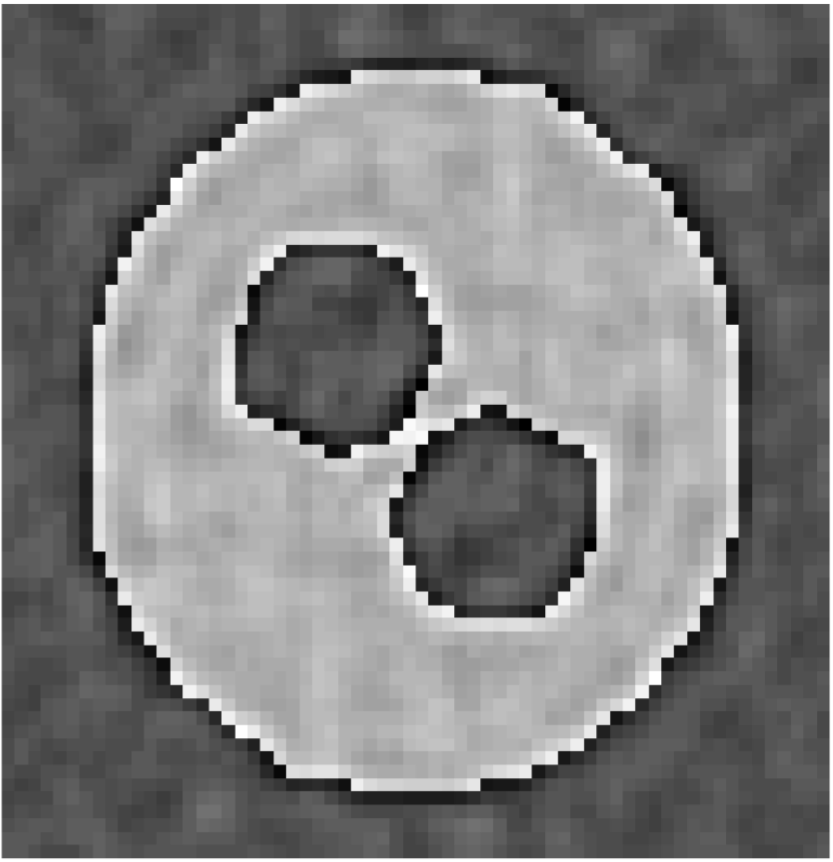}
        \caption{$\delta = 2.5 $}
        \label{fig:}
    \end{subfigure}
 \caption{We show the reconstructions obtained solving \eqref{eq:zerojoint} for different values of $\delta$. For $\delta=0$ we get the zero-filling solution. For small $\delta$ we expect the solution to  be similar to the zero-filling reconstruction. For $\delta=1$ we see the effect of the joint term on the reconstruction. The solution presents the same noise artefacts but having in addition very sharp  boundaries. Finally, for very large $\delta$ we still have enhanced boundaries but we also amplify the noise. }
\label{fig:deltas}
\end{figure}

In the model proposed in \eqref{eq:joint}, the parameters that we need to choose are $\alpha$, $\beta$ and $\delta$. In this section we discuss a rule to choose them depending on the desired results. Some examples will clarify these empirical choices. 
\begin{itemize}
\item $\alpha$ balances the total variation regularization term in the reconstruction for the magnitude. The higher the $\alpha$, the more piecewise constant the reconstruction will be. See \autoref{fig:alphas}.

\item $\beta$ defines the scale of the objects that will be detected in the segmentation. Smaller values of $\beta$ will allow for smaller objects. See  \autoref{fig:betas}.

\item $\delta$ is the parameter linking the reconstruction and the segmentation. To better illustrate its role, let us consider   a zero-filling like reconstruction and segmentation:
\begin{equation}
 \imath_{\mathcal{SF} \cdot = f}(u)+ \delta \sum_i\in\Omega  \sum_{j=1}^{\ell} v_{ij} (c_j - u_i )^2 dx + \beta \|\nabla v\| \rightarrow \min_{u,v}
\label{eq:zerojoint}
\end{equation}
where $ \imath(u) = \begin{cases}
      +\infty , & \text{if}\ \mathcal{SF}u \neq f \\
      0, &  \text{if}\ \mathcal{SF}u = f
    \end{cases} $. This problem is solving the zero-filled reconstruction and segmentation jointly. For $\delta=0$, the reconstruction is the zero-filling solution. In \autoref{fig:deltas} we can see the impact of the segmentation term on the reconstruction for increasing values of $\delta$. We can see that for very small$\delta$ the result is close to the zero-filling solution. For $\delta=1$ the noise from the model is present as expected but in addition the boundaries are enhanced. For large $\delta$ the boundaries are still very pronounced and the noise is also amplified. 
\end{itemize}
\subsection{Comparison with another joint approach}
\begin{figure}
\centering
\begin{subfigure}{0.2\textwidth}
 \includegraphics[width=\textwidth]{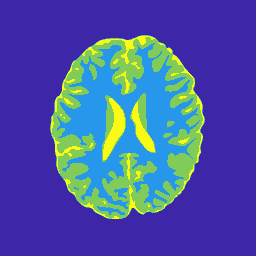}
        \caption{Groundtruth \newline}
        \label{fig:}
    \end{subfigure}
~
\begin{subfigure}{0.2\textwidth}
 \includegraphics[width=\textwidth]{brain_jointbregseg1.png}
        \caption{Joint segmentation\newline} 
        \label{fig:}
    \end{subfigure}
    ~
\begin{subfigure}{0.2\textwidth}
 \includegraphics[width=\textwidth]{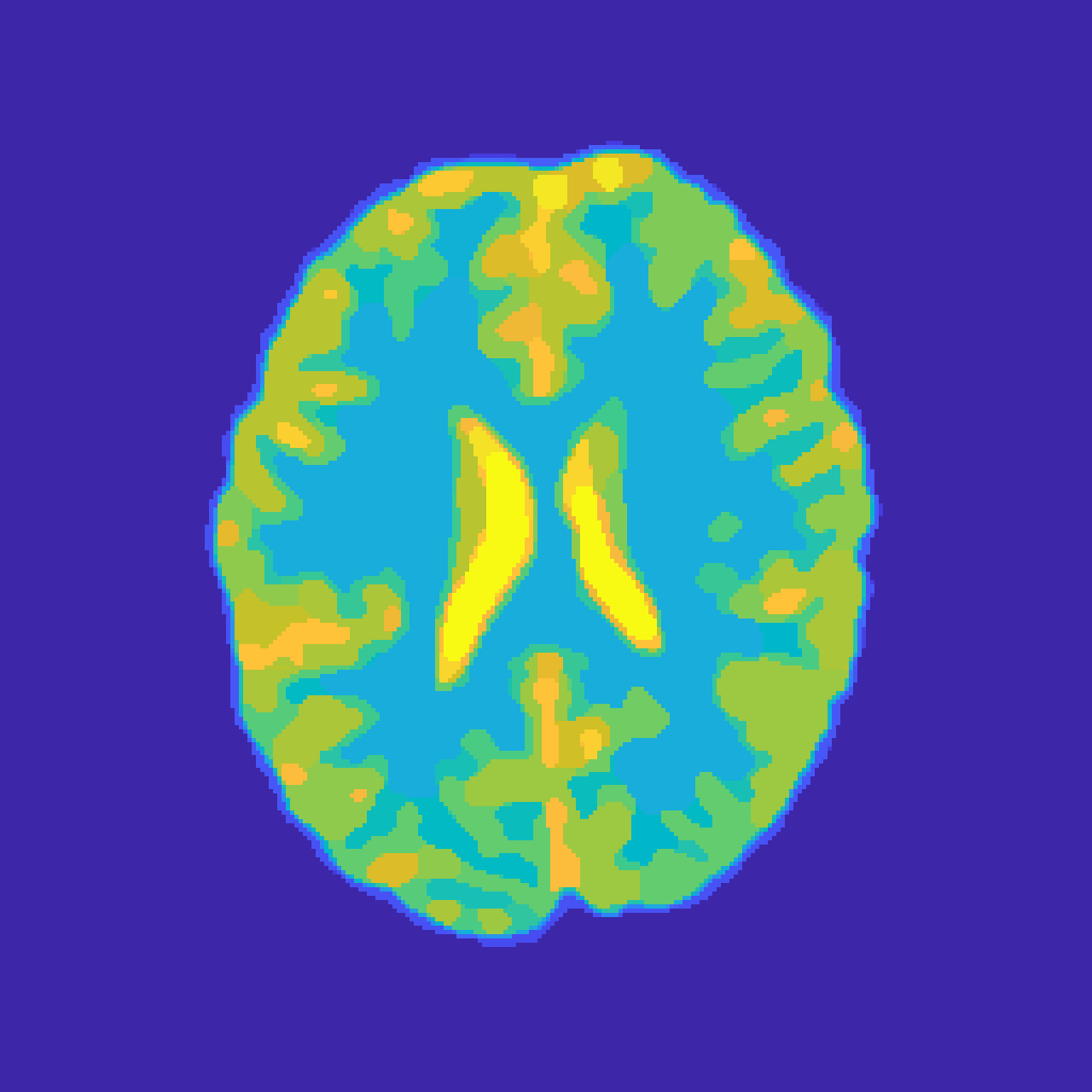}
        \caption{Potts model, $\gamma=0.01$} 
        \label{fig:}
    \end{subfigure}
        ~
\begin{subfigure}{0.2\textwidth}
 \includegraphics[width=\textwidth]{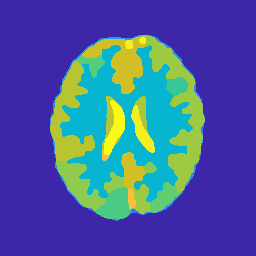}
        \caption{Potts model, $\gamma=0.05$} 
        \label{fig:}
    \end{subfigure}

    \begin{subfigure}{0.2\textwidth}
 \includegraphics[width=\textwidth]{bubbleNew_gt.png}
        \caption{Groundtruth \newline }
        \label{fig:}
    \end{subfigure}
~
\begin{subfigure}{0.2\textwidth}
 \includegraphics[width=\textwidth]{bubbleNew1_jointbregseg1.png}
        \caption{Joint segmentation \newline}
        \label{fig:}
    \end{subfigure}
    ~
\begin{subfigure}{0.2\textwidth}
 \includegraphics[width=\textwidth]{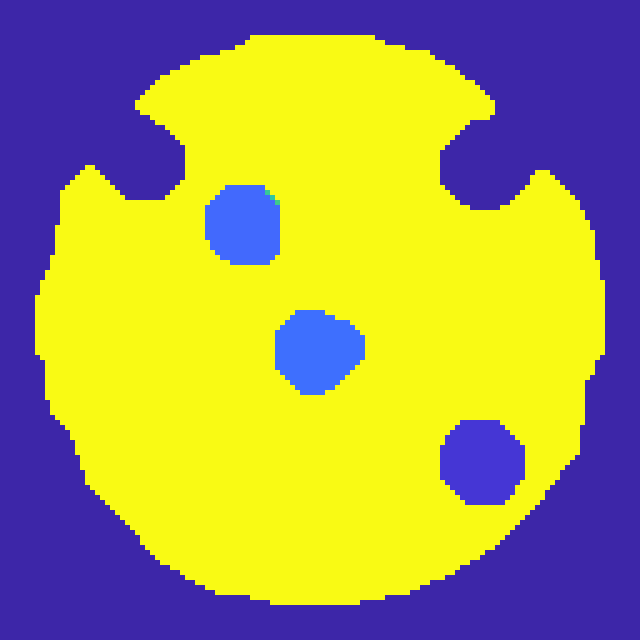}
        \caption{Potts model,  $\gamma=0.75$}
        \label{fig:}
    \end{subfigure}
        ~
\begin{subfigure}{0.2\textwidth}
 \includegraphics[width=\textwidth]{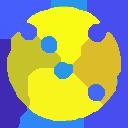}
        \caption{Potts model,  $\gamma=0.5$ \newline}
        \label{fig:}
    \end{subfigure}
        
    \begin{subfigure}{0.2\textwidth}
    \includegraphics[width=\textwidth]{tum_zero_seg.png}
        \caption{Segmentation from zero-filled reconstruction}
        \label{fig:}
    \end{subfigure}
~
\begin{subfigure}{0.2\textwidth}
 \includegraphics[width=\textwidth, height=3.21cm]{tum_joint_seg.png}
        \caption{Joint segmentation \newline\newline}
        \label{fig:}
    \end{subfigure}
    ~
\begin{subfigure}{0.2\textwidth}
 \includegraphics[width=\textwidth]{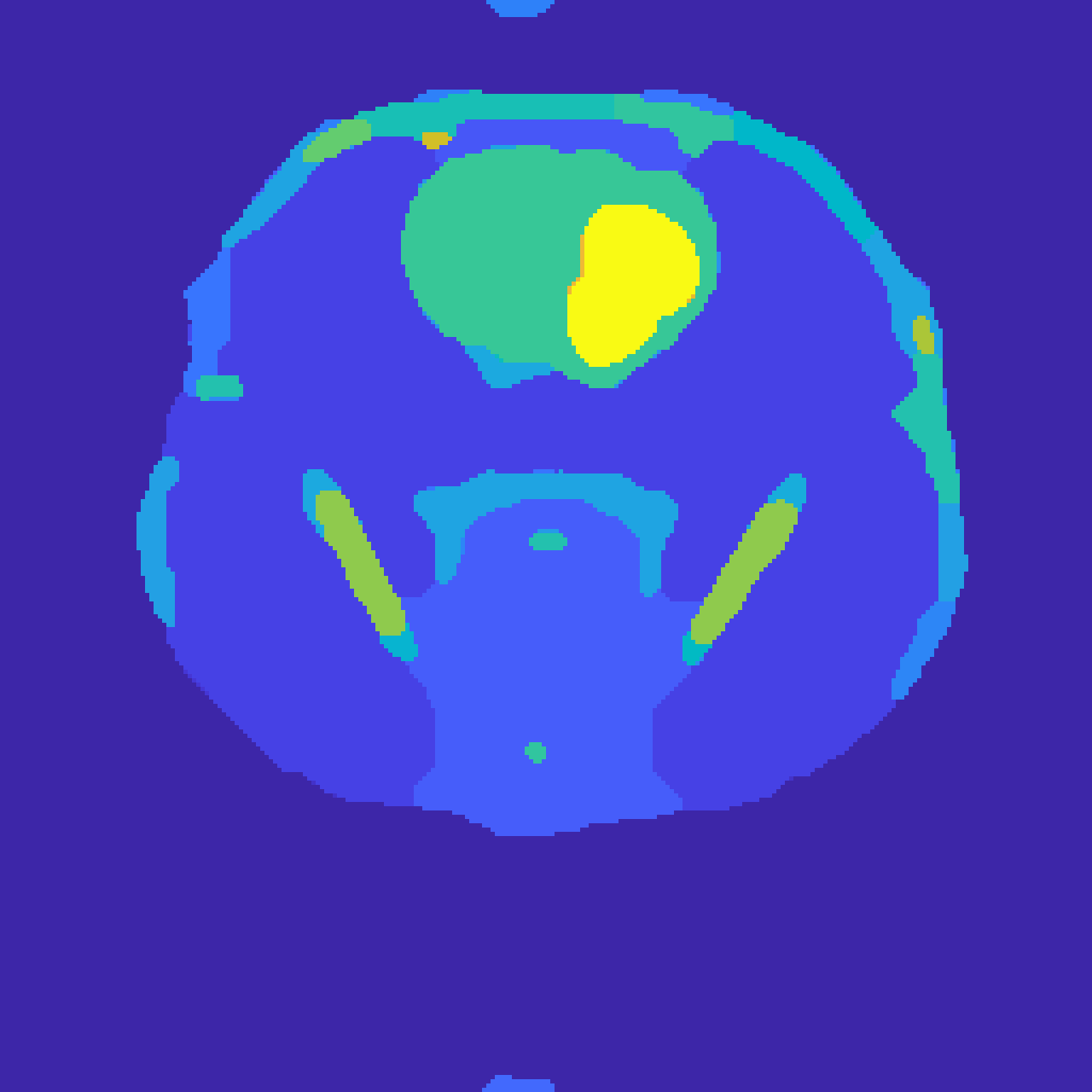}
        \caption{Potts model, $\gamma=0.05$\newline }
        \label{fig:}
    \end{subfigure}
        ~
\begin{subfigure}{0.2\textwidth}
 \includegraphics[width=\textwidth]{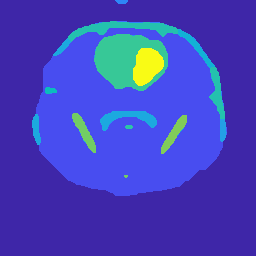}
        \caption{Potts model, $\gamma=0.1$ \newline\newline}
        \label{fig:}
    \end{subfigure}
    \caption{Comparison of our joint approach with the Potts model. Noise level and undersampling rate are described in \autoref{fig:brain_whole}, \ref{fig:bubblyflow} and \ref{fig:tum}. The results are presented for three different examples and for two different choices of the regularisation parameter $\gamma$. We can see that the Potts model tends to overestimate the number of regions to segment.}
    \label{fig:comparisonPotts}
\end{figure}

We present a comparison of our joint model with another non-convex method, namely the Potts model approach by \cite{storath}, described in \autoref{subsec:comparison}. The major advantage of the joint reconstruction and segmentation using the Potts model is that it does not require to select explicitely the number of regions to segment, although this depends on the choice of the regularisation parameter. However, by definition, it only produces a piecewise constant image, therefore a segmentation, and not a reconstruction. This is useful in some applications where one is only interested in the segmentation. In  contrast, our model produces both reconstruction and segmentation. In \autoref{fig:comparisonPotts}, we show the results for some examples. Note that because the results of the Potts model are in the range of the groundtruth image, while our segmentation are in label space, we can not directly use the RSE as before, or common metrics that compare actual intensities such as PSNR and structure similarity index measure (SSIM), for comparison. For example, for some tissue in class 1, to label it class 2 is as wrong as to label it class 3. However in this case, the SSIM and PSNR will favour the label class 2. \\
We therefore focus on a visual assessment and show the results of the Potts model for two different choices of the regularisation parameter $\gamma$. We recall that the proposed model requires to determine the number of classes in advance, while the model for comparison estimates the number of regions but this depends on the choice of the regularisation parameter. In the top row, we can see that the Potts model, although it retrieves the shape of the main structures for the brain phantom example, it overestimates the number of classes. By increasing the parameter $\gamma$, this issue is not resolved as it assigns different intensities to objects of the same category. In contrast, our approach is able to identify the desired classes as in the groundtruth. For the bubble case (middle row), we can see that our method works better and our segmentation is more accurate, while the Potts model fails to capture shape details (e.g., outer circle is distorted) and again overestimates the number of regions. We can also see that, when slightly decreasing $\gamma$, the Potts model is very sensitive to artefacts. For the real MR data (bottom row), we see that both methods identify the tumour quite well. Because we were only interested in identifying three classes as tumour, brain and background, we do not segment the outer region (rat's head), captured insted by the Potts model. However, the Potts model only produces the segmentation, while our method, as shown in \autoref{fig:tum}, also produces an enhanced reconstruction with sharp edges.

\section{Conclusion}
\label{sec:conclusion}
In this paper, we have investigated a novel mathametical approach to perform simultaneously reconstruction and segmentation from undersampled MRI data. Our motivation was to include in the reconstruction prior knowledge of the objects we are interested in. By interconnecting the reconstruction and the segmentation terms, we can achieve sharper reconstructions and more accurate segmentations. We derived a variational model based on Bregman iteration and we have verified its convergence properties. With our approach we show that by solving the more complicated joint model, we are able to improve both reconstruction and segmentation compared to the traditional sequential approach. This suggests that with the joint model it is possible to reduce error propagations that occur in sequential analysis, when the segmentation is separate and posterior to the reconstruction. \\
\indent We have tested our method for two different application, which are bubbly flow and cancer imaging. In both cases, the reconstructions are sharper and finer structures are detected. Additionally, the segmentations also benefit from the improvement in the reconstructions. We have found that the joint model outperforms the sequential approach by exploiting prior information on the objects that we want to segment. In addition, we also show that our method performs better than the well-known Potts model. We also presented a discussion on the parameter choice rule that offer some insight on how to tune the parameters according to the desired result. It is interesting to notice that, with our model, we are able to control the segmentation effect on the reconstruction. Furthermore, when the final analysis of the MR image is indeed the segmentation, it is possible to bias the reconstruction towards the piecewise constant solution, yet preserving finer details in the structure. \\
\indent In our set-up, we have specified the intensity constants characteristic of the region of interests, which were known a priori for our applications. However, it is possible to also include the optimisation with respect to $c_j$ in our joint model, where the same convergence guarantees hold (see Remark 2). Nevertheless, one limitation of the model is the need to specify the number of regions to be segmented.\\
\indent In our future research, we would like to study the extension of this model for the bubbly flow to the reconstruction of the magnitude image as well as the phase image. The goal is not only to extract the structure of the bubble, but also to estimate velocity information, which is encoded in the phase image. As the problem is non-convex in the joint argument but also non-convex with respect to the phase, we need to derive a different convergence analysis. 

\section*{Acknowledgements}
VC acknowledges the financial support of the Cambridge Cancer Centre. MB acknowledges the Leverhulme Trust Early Career Fellowship "Learning from mistakes: a supervised feedback-loop for imaging applications", the Isaac Newton Trust and the Cantab Capital Institute for the Mathematics of Information. MJE and CBS acknowledge support from Leverhulme Trust project "Breaking the non-convexity barrier",
EPSRC grant "EP/M00483X/1", EPSRC centre "EP/N014588/1", the Cantab Capital Institute for the Mathematics of Information,
and from CHiPS and NoMADS (Horizon 2020 RISE project grant). Moreover, CBS is thankful for support by the Alan Turing Institute. We would like to thank Kevin Brindle and Alan Wright for the acquisition of the data in \autoref{fig:tum}.

\section*{References}
  \addcontentsline{toc}{section}{Bibliography}
  \bibliographystyle{unsrt}
  \bibliography{alternating} 

\begin{thebibliography}{10}

\bibitem{Candes2006}
Emmanuel~J. Cand{\`e}s, Justin Romberg, and Terence Tao.
\newblock Robust uncertainty principles: Exact signal reconstruction from
  highly incomplete frequency information.
\newblock {\em IEEE Transactions on information theory}, 52(2):489--509, 2006.

\bibitem{Donoho2006}
David~L. Donoho.
\newblock Compressed sensing.
\newblock {\em IEEE Transactions on information theory}, 52(4):1289--1306,
  2006.

\bibitem{LustigDonohoPauly2007}
Michael Lustig, David Donoho, and John~M. Pauly.
\newblock Sparse {MRI}: The application of compressed sensing for rapid {MR}
  imaging.
\newblock {\em Magnetic resonance in medicine}, 58(6):1182--1195, 2007.

\bibitem{Macovski1996}
Albert Macovski.
\newblock Noise in mri.
\newblock {\em Magnetic Resonance in Medicine}, 36(3):494--497, 1996.

\bibitem{Ehrhardt2016mri}
Matthias~J. Ehrhardt and Marta~M. Betcke.
\newblock Multi-contrast {MRI} reconstruction with structure-guided total
  variation.
\newblock {\em SIAM Journal on Imaging Sciences}, 9(3):1084--1106, 2016.

\bibitem{BenningGladdenHollandEtAl2014}
Martin Benning, Lynn Gladden, Daniel Holland, Carola-Bibiane Sch{\"o}nlieb, and
  Tuomo Valkonen.
\newblock Phase reconstruction from velocity-encoded {MRI} measurements--a
  survey of sparsity-promoting variational approaches.
\newblock {\em Journal of Magnetic Resonance}, 238:26--43, 2014.

\bibitem{Meyer2001}
Yves Meyer.
\newblock {\em Oscillating patterns in image processing and nonlinear evolution
  equations: the fifteenth Dean Jacqueline B. Lewis memorial lectures},
  volume~22.
\newblock American Mathematical Soc., 2001.

\bibitem{StrongChan2003}
David Strong and Tony Chan.
\newblock Edge-preserving and scale-dependent properties of total variation
  regularization.
\newblock {\em Inverse problems}, 19(6):S165, 2003.

\bibitem{groundstates}
Martin Benning and Martin Burger.
\newblock {G}round states and singular vectors of convex variational
  regularization methods.
\newblock {\em Methods and Applications of Analysis}, 2012.

\bibitem{Osher2005}
Stanley Osher, Martin Burger, Donald Goldfarb, Jinjun Xu, and Wotao Yin.
\newblock An iterative regularization method for total variation-based image
  restoration.
\newblock {\em Multiscale Modeling \& Simulation}, 4(2):460--489, 2005.

\bibitem{breg}
Lev~M. Bregman.
\newblock {T}he relaxation method of finding the common points of convex sets
  and its application to the solution of problems in convex optimization.
\newblock {\em USSR Computational Mathematics and Mathematical Physics}, 1967.

\bibitem{Kiwiel1997}
Krzysztof~C. Kiwiel.
\newblock Proximal minimization methods with generalized bregman functions.
\newblock {\em SIAM journal on control and optimization}, 35(4):1142--1168,
  1997.

\bibitem{ROF}
Leonid~I. Rudin, Stanley Osher, and Emad Fatemi.
\newblock {N}onlinear total variation based noise removal algorithms.
\newblock {\em Physica D, vol. 60, pp. 259-268}, 1992.

\bibitem{Morozov1966}
Vladimir~Alekseevich Morozov.
\newblock On the solution of functional equations by the method of
  regularization.
\newblock In {\em Soviet Math. Dokl}, volume~7, pages 414--417, 1966.

\bibitem{Mumford1989}
David Mumford and Jayant Shah.
\newblock Optimal approximations by piecewise smooth functions and associated
  variational problems.
\newblock {\em Communications on pure and applied mathematics}, 42(5):577--685,
  1989.

\bibitem{Ambrosio1990}
Luigi Ambrosio and Vincenzo~Maria Tortorelli.
\newblock Approximation of functional depending on jumps by elliptic functional
  via t-convergence.
\newblock {\em Communications on Pure and Applied Mathematics},
  43(8):999--1036, 1990.

\bibitem{Chambolle1995}
Antonin Chambolle.
\newblock Image segmentation by variational methods: Mumford and shah
  functional and the discrete approximations.
\newblock {\em SIAM Journal on Applied Mathematics}, 55(3):827--863, 1995.

\bibitem{Boysen2009}
Leif Boysen, Angela Kempe, Volkmar Liebscher, Axel Munk, Olaf Wittich, et~al.
\newblock Consistencies and rates of convergence of jump-penalized least
  squares estimators.
\newblock {\em The Annals of Statistics}, 37(1):157--183, 2009.

\bibitem{Pock2009}
Thomas Pock, Antonin Chambolle, Daniel Cremers, and Horst Bischof.
\newblock A convex relaxation approach for computing minimal partitions.
\newblock In {\em 2009 IEEE Conference on Computer Vision and Pattern
  Recognition}, pages 810--817, June 2009.

\bibitem{Chant2001}
Tony~F. Chan and Luminita~A. Vese.
\newblock Active contours without edges.
\newblock {\em IEEE Transactions on Image Processing}, 10(2):266--277, 2001.

\bibitem{Ambrosio2000}
Luigi Ambrosio, Nicola Fusco, and Diego Pallara.
\newblock {\em Functions of bounded variation and free discontinuity problems},
  volume 254.
\newblock Clarendon Press Oxford, 2000.

\bibitem{Chan2006}
Tony~F. Chan, Selim Esedoglu, and Mila Nikolova.
\newblock Algorithms for finding global minimizers of image segmentation and
  denoising models.
\newblock {\em SIAM journal on applied mathematics}, 66(5):1632--1648, 2006.

\bibitem{Lellmann2013}
Jan Lellmann, Bj{\"o}rn Lellmann, Florian Widmann, and Christoph Schn{\"o}rr.
\newblock Discrete and continuous models for partitioning problems.
\newblock {\em International journal of computer vision}, 104(3):241--269,
  2013.

\bibitem{Lellmann2009}
Jan Lellmann, J{\"o}rg Kappes, Jing Yuan, Florian Becker, and Christoph
  Schn{\"o}rr.
\newblock Convex multi-class image labeling by simplex-constrained total
  variation.
\newblock In {\em International conference on scale space and variational
  methods in computer vision}, pages 150--162. Springer, 2009.

\bibitem{Zeune2017}
Leonie Zeune, Guus van Dalum, Leon W. M.~M. Terstappen, Stephan~A. van Gils,
  and Christoph Brune.
\newblock Multiscale segmentation via bregman distances and nonlinear spectral
  analysis.
\newblock {\em SIAM journal on imaging sciences}, 10(1):111--146, 2017.

\bibitem{ramlau}
Ronny Ramlau and Wolfgang Ring.
\newblock {A} {M}umford-{S}hah level-set approach for the inversion and
  segmentation of {X}-ray tomography data.
\newblock {\em Journal of Computational Physics, 221, 2, 539 - 557}, 2007.

\bibitem{klann}
Ronny~Ramlau Esther~Klann and Wolfgang Ring.
\newblock {A} {M}umford-{S}hah level-set approach for the inversion and
  segmentation of {SPECT}/{CT} data.
\newblock {\em Inverse Problems and Imaging, 5, 1, 137 - 166}, 2011.

\bibitem{Klann1}
Esther Klann.
\newblock {A} {M}umford-{S}hah-like method for limited data tomography with an
  application to electron tomography.
\newblock {\em SIAM J. Imaging Sci., 4(4), 1029–1048}, 2011.

\bibitem{0266-5611-32-10-104002}
Martin Burger, Carolin Rossmanith, and Xiaoqun Zhang.
\newblock {S}imultaneous reconstruction and segmentation for dynamic {SPECT}
  imaging.
\newblock {\em Inverse Problems}, 32(10):104002, 2016.

\bibitem{Lauze2017}
Fran{\c{c}}ois Lauze, Yvain Qu{\'e}au, and Esben Plenge.
\newblock Simultaneous reconstruction and segmentation of {CT} scans with
  shadowed data.
\newblock In {\em International Conference on Scale Space and Variational
  Methods in Computer Vision}, pages 308--319. Springer, 2017.

\bibitem{van}
Dominique Van~de Sompel and Michael Brady.
\newblock {S}imultaneous reconstruction and segmentation algorithm for positron
  emission tomography and transmission tomography.
\newblock {\em 5th IEEE International Symposium on Biomedical Imaging: From
  Nano to Macro, 1035-1038}, 2008.

\bibitem{yiqiu}
Mikhail Romanov, Anders~Bjorholm Dahl, Yiqiu Dong, and Per~Christian Hansen.
\newblock {S}imultaneous tomographic reconstruction and segmentation with class
  priors.
\newblock {\em Inverse Problems in Science and Engineering, 24,8, 1432-1453},
  2015.

\bibitem{storath}
Martin Storath, Andreas Weinmann, J{\"u}rgen Frikel, and Michael Unser.
\newblock {J}oint image reconstruction and segmentation using the {P}otts
  model.
\newblock {\em Inverse Problems 31 (2015) 025003 (29pp)}, 2015.

\bibitem{caballero}
Jose Caballero, Wenjia Bai, Anthony~N. Price, Daniel Rueckert, , and Joseph~V.
  Hajnal.
\newblock {A}pplication-driven {MRI}: {J}oint reconstruction and segmentation
  from undersampled {MRI} data.
\newblock {\em Med Image Comput Comput Assist Interv;17(Pt 1):106-13.}, 2014.

\bibitem{pd}
Antonin Chambolle and Thomas Pock.
\newblock {A} first-order primal-dual algorithm for convex problems with
  applications to imaging.
\newblock {\em J Math Imaging Vis 40: 120. doi:10.1007/s10851-010-0251-1},
  2011.

\bibitem{ChambollePock2016}
Antonin Chambolle and Thomas Pock.
\newblock An introduction to continuous optimization for imaging.
\newblock {\em Acta Numerica}, 25:161--319, 2016.

\bibitem{EsserZhangChan2010}
Ernie Esser, Xiaoqun Zhang, and Tony~F. Chan.
\newblock A general framework for a class of first order primal-dual algorithms
  for convex optimization in imaging science.
\newblock {\em SIAM Journal on Imaging Sciences}, 3(4):1015--1046, 2010.

\bibitem{PockCremersBischofEtAl2009}
Thomas Pock, Daniel Cremers, Horst Bischof, and Antonin Chambolle.
\newblock An algorithm for minimizing the {M}umford-{S}hah functional.
\newblock In {\em Computer Vision, 2009 IEEE 12th International Conference on},
  pages 1133--1140. IEEE, 2009.

\bibitem{ChambolleCremersPock2012}
Antonin Chambolle, Daniel Cremers, and Thomas Pock.
\newblock A convex approach to minimal partitions.
\newblock {\em SIAM Journal on Imaging Sciences}, 5(4):1113--1158, 2012.

\bibitem{Lojasiewicz1963}
Stanislaw Lojasiewicz.
\newblock Une propri{\'e}t{\'e} topologique des sous-ensembles analytiques
  r{\'e}els.
\newblock {\em Les {\'e}quations aux d{\'e}riv{\'e}es partielles}, 117:87--89,
  1963.

\bibitem{BenningBetckeEhrhardtEtAl2017}
Martin Benning, Marta~M. Betcke, Matthias~J. Ehrhardt, and Carola-Bibiane
  Sch{\"o}nlieb.
\newblock Choose your path wisely: gradient descent in a {B}regman distance
  framework.
\newblock {\em arXiv preprint arXiv:1712.04045}, 2017.

\bibitem{PALM}
Jerome Bolte, Shoham Sabach, and Marc Teboulle.
\newblock Proximal alternating minimization for nonconvex and nonsmooth
  problems.
\newblock {\em Mathematical Programming, 146(1-2):459-494}, 2014.

\bibitem{Chalmers1994}
Jeffrey~J. Chalmers.
\newblock Cells and bubbles in sparged bioreactors.
\newblock In {\em Cell Culture Engineering IV}, pages 311--320. Springer, 1994.

\bibitem{Deckwer1992}
Wolf-Dieter Deckwer and Robert~W. Field.
\newblock {\em Bubble column reactors}, volume 200.
\newblock Wiley New York, 1992.

\bibitem{Holland2012}
Daniel~J. Holland, Andrew Blake, Alexander~B. Tayler, Andrew~J. Sederman, and
  Lynn~F. Gladden.
\newblock Bubble size measurement using bayesian magnetic resonance.
\newblock {\em Chemical engineering science}, 84:735--745, 2012.

\bibitem{Tayler2012}
Alexander~B. Tayler, Daniel~J. Holland, Andrew~J. Sederman, and Lynn~F.
  Gladden.
\newblock Applications of ultra-fast {MRI} to high voidage bubbly flow:
  measurement of bubble size distributions, interfacial area and hydrodynamics.
\newblock {\em Chemical engineering science}, 71:468--483, 2012.

\bibitem{HollandMalioutovBlakeEtAl2010}
Daniel~J. Holland, Dmitry~M. Malioutov, Andrew Blake, Andrew~J. Sederman, and
  L.~F. Gladden.
\newblock Reducing data acquisition times in phase-encoded velocity imaging
  using compressed sensing.
\newblock {\em Journal of magnetic resonance}, 203(2):236--246, 2010.

\bibitem{TaylerHollandSedermanEtAl2012}
Alexander~B. Tayler, Daniel~J. Holland, Andrew~J. Sederman, and Lynn~F.
  Gladden.
\newblock Exploring the origins of turbulence in multiphase flow using
  compressed sensing {MRI}.
\newblock {\em Physical review letters}, 108(26):264505, 2012.

\bibitem{Rodrigues2014}
Tiago~B. Rodrigues, Eva~M. Serrao, Brett W.~C. Kennedy, De-En Hu, Mikko~I.
  Kettunen, and Kevin~M. Brindle.
\newblock Magnetic resonance imaging of tumor glycolysis using hyperpolarized
  13 c-labeled glucose.
\newblock {\em Nature medicine}, 20(1):93, 2014.

\end{thebibliography}
\clearpage
\appendix

\section{Numerical results on phantoms}
\label{phantoms}
%
%
%
%
%

\begin{figure}[h]
\centering
\begin{subfigure}{0.23\textwidth}
\centering
 \includegraphics[width=\textwidth]{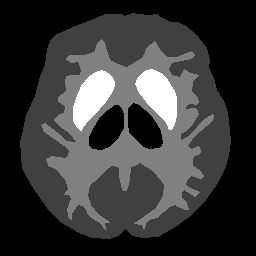}
        \caption{Groundtruth \newline \newline}
        \label{fig:gtbrain}
    \end{subfigure}
~
\begin{subfigure}{0.23\textwidth}
\centering
 \includegraphics[width=\textwidth]{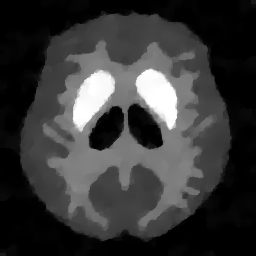}
        \caption{TV reconstruction, $\alpha=0.15$, RRE=0.0305, PSNR=27.44}
        \label{fig:brain256TV}
    \end{subfigure}
    ~
\begin{subfigure}{0.23\textwidth}
\centering
 \includegraphics[width=\textwidth]{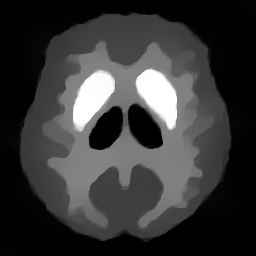}
        \caption{Bregman reconstruction, $\alpha=1.1$, RRE=0.0427, PSNR=27.21}
        \label{fig:brain256Breg}
    \end{subfigure}
~
    \begin{subfigure}{0.23\textwidth}
    \centering
 \includegraphics[width=\textwidth]{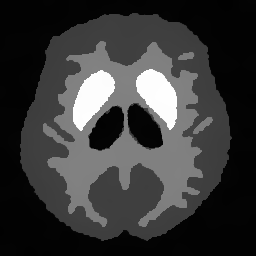}
        \caption{Joint reconstruction, $\alpha=0.8$,  RRE=0.0262, PSNR=28.27 }
        \label{fig:brain256Joint}
    \end{subfigure}

\begin{subfigure}{0.23\textwidth}
\centering
 \includegraphics[width=\textwidth]{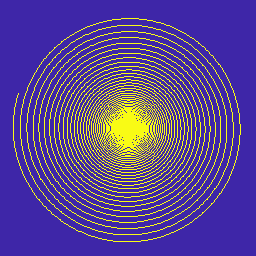}
        \caption{Sampling matrix, 15\% \newline}
        \label{fig:sampling15}
    \end{subfigure}
    ~
            \begin{subfigure}{0.23\textwidth}
        \centering
    \includegraphics[width=\textwidth]{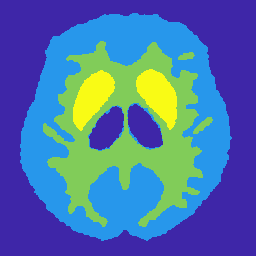}
        \caption{Segmentation, $\beta=0.001$  RSE=0.0219}
        \label{fig:brain256TVseg}
    \end{subfigure}
        ~
\begin{subfigure}{0.23\textwidth}
\centering
 \includegraphics[width=\textwidth]{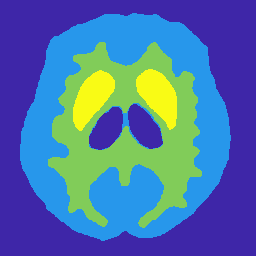}
        \caption{Bregman segmentation, $\beta=0.001$ RSE=0.0399}
        \label{fig:brain256Bregseg}
    \end{subfigure}
    ~
    \begin{subfigure}{0.23\textwidth}
    \centering
 \includegraphics[width=\textwidth]{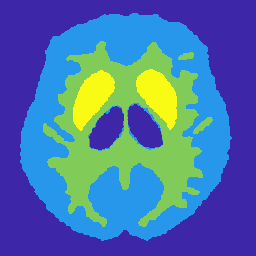}
        \caption{Joint segmentation, $\beta=0.001$, $\delta=2$, RSE=0.0219}
        \label{fig:brain256Jointseg}
    \end{subfigure} 
    \caption{This example shows clearly the effect of the parameter $\delta$ in the joint model. The segmentation is easy to achieve and we do not see a significant improvement in joint segmentation compared to the TV sequential segmentation, but there is a small gain compared to the sequential Bregman segmantation. However, the joint reconstruction results improved thanks to the parameter $\delta$ which biases the reconstruction to be closer to the segmentation. }
    \label{fig:brain256}
\end{figure}

\begin{figure}[h]
\centering
\begin{subfigure}{0.23\textwidth}
\centering
 \includegraphics[width=\textwidth]{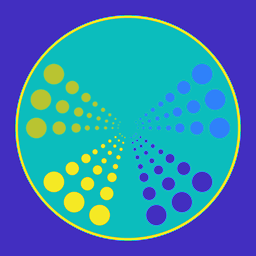}
        \caption{Groundtruth \newline \newline}
        \label{fig:gtcircle}
    \end{subfigure}
~
\begin{subfigure}{0.23\textwidth}
\centering
 \includegraphics[width=\textwidth]{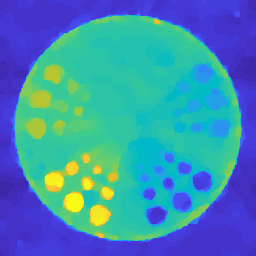}
        \caption{TV reconstruction, $\alpha=0.3$, RRE=0.0578, PSNR=21.43}
        \label{fig:circleTV}
    \end{subfigure}
    ~
\begin{subfigure}{0.23\textwidth}
\centering
 \includegraphics[width=\textwidth]{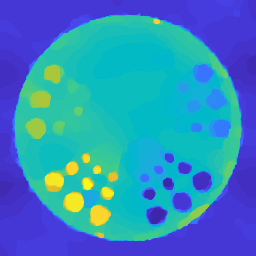}
        \caption{Bregman reconstruction, $\alpha=1.5$, RRE= 0.1307, PSNR=21.49}
        \label{fig:circleBreg}
    \end{subfigure}
~
    \begin{subfigure}{0.23\textwidth}
    \centering
 \includegraphics[width=\textwidth]{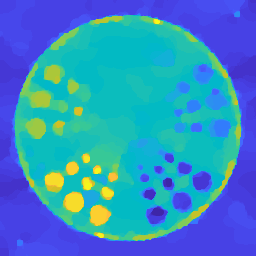}
        \caption{Joint reconstruction, $\alpha=1.5$ RRE=0.0713, PSNR=21.87}
        \label{fig:circleJoint}
    \end{subfigure}

\begin{subfigure}{0.23\textwidth}
\centering
 \includegraphics[width=\textwidth]{presentation/images/map15_col.png}
        \caption{Sampling matrix, 15\% \newline}
        \label{fig:sampling15}
    \end{subfigure}
    ~
            \begin{subfigure}{0.23\textwidth}
        \centering
    \includegraphics[width=\textwidth]{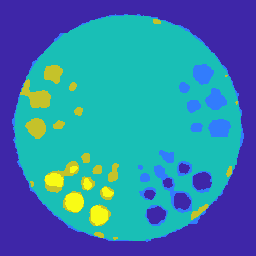}
        \caption{Segmentation, $\beta=0.001$, RSE=0.096}
        \label{fig:circleTVseg}
    \end{subfigure}
        ~
\begin{subfigure}{0.23\textwidth}
\centering
 \includegraphics[width=\textwidth]{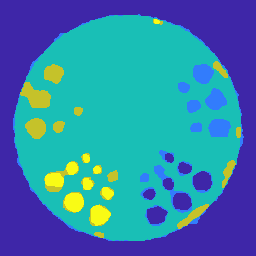}
        \caption{Bregman segmentation, $\beta=0.001$, RSE=0.121}
        \label{fig:circleBregseg}
    \end{subfigure}
    ~
    \begin{subfigure}{0.23\textwidth}
    \centering
 \includegraphics[width=\textwidth]{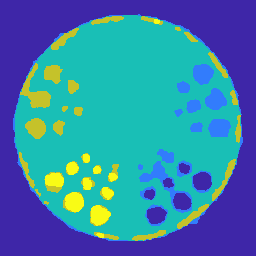}
        \caption{Joint segmentation, $\beta=0.001$, $\delta=0.1$, RSE=0.091}
        \label{fig:circleJointseg}
    \end{subfigure} 
    \caption{In this example, we can see that the reconstructions are quite similar. However in the joint reconstruction, the outer yellow circle, which is completely ignored by the sequential reconstructions, is partially detected. This is also the case for the joint segmenation.}
    \label{fig:circleflow}
\end{figure}


\begin{figure}[h]
\centering
\begin{subfigure}{0.23\textwidth}
\centering
 \includegraphics[width=\textwidth]{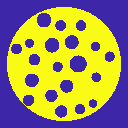}
        \caption{Groundtruth \newline \newline}
        \label{fig:gt1}
    \end{subfigure}
~
\begin{subfigure}{0.23\textwidth}
\centering
 \includegraphics[width=\textwidth]{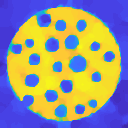}
        \caption{TV reconstruction, $\alpha=0.15$, RRE=0.074, PSNR=17.15 }
        \label{fig:bubblyTV1}
    \end{subfigure}
    ~
\begin{subfigure}{0.23\textwidth}
\centering
 \includegraphics[width=\textwidth]{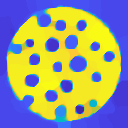}
        \caption{Bregman reconstruction, $\alpha=2$,  RRE=0.071, PSNR=17.65}
        \label{fig:bubblyBreg1}
    \end{subfigure}
~
    \begin{subfigure}{0.23\textwidth}
    \centering
 \includegraphics[width=\textwidth]{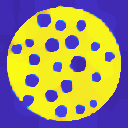}
        \caption{Joint reconstruction, $\alpha=0.8$, RRE=0.047, PSNR=19.015 } 
        \label{fig:bubblyJoint1}
    \end{subfigure}

\begin{subfigure}{0.23\textwidth}
\centering
 \includegraphics[width=\textwidth]{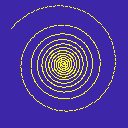}
        \caption{Sampling matrix, 8\% \newline
        }
        \label{fig:sampling1}
    \end{subfigure}
    ~
            \begin{subfigure}{0.23\textwidth}
        \centering
    \includegraphics[width=\textwidth]{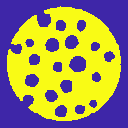}
        \caption{Segmentation, $\beta=0.01$, RSE=0.016 }
        \label{fig:bubblyTVseg1}
    \end{subfigure}
        ~
\begin{subfigure}{0.23\textwidth}
\centering
 \includegraphics[width=\textwidth]{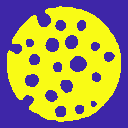}
        \caption{Bregman segmentation, $\beta=0.01$, RSE=0.022}
        \label{fig:bubbleBregseg1}
    \end{subfigure}
    ~
    \begin{subfigure}{0.23\textwidth}
    \centering
 \includegraphics[width=\textwidth]{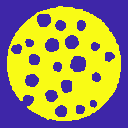}
        \caption{Joint segmentation, $\beta=0.01$, $\delta=1$, RSE=0.014}
        \label{fig:bubblyJointseg1}
    \end{subfigure} 
    \caption{In this example for the bubbly flow, we can see clearly an improvement for both joint reconstruction and joint segmentation. The contrast in the joint reconstruction is better recovered and the segmentation is more accurate, especially for the bubbles close to the edge of the pipe. The joint method results particularly useful for the bubbly flow application.}
    \label{fig:bubblyflow1}
\end{figure}
\end{document}